\documentclass{amsart}
\usepackage{eurosym}
\usepackage{amssymb,amsmath,amsfonts}
\usepackage{amsmath}
\usepackage{tikz}
\usepackage[all]{xy}
\usepackage[unicode=true]{hyperref}

\newtheorem{theorem}{Theorem}[section]
\newtheorem{corollary}[theorem]{Corollary}

\newtheorem{definition}[theorem]{Definition}
\newtheorem{example}[theorem]{Example}

\newtheorem{lemma}[theorem]{Lemma}
\newtheorem{remark}[theorem]{Remark}
\newtheorem{proposition}[theorem]{Proposition}

\newcommand{\cref}[1]{eq. (\ref{#1})}

\begin{document}
\title[Matrix invertible extensions]{Matrix invertible
extensions over commutative rings. Part II: determinant liftability}
\author{Grigore C\u{a}lug\u{a}reanu, Horia F.\ Pop, Adrian Vasiu}

\date{Accepted for publication in final form in Linear Algebra Appl. on July 5, 2025.}

\begin{abstract}
A unimodular $2\times 2$ matrix $A$ with entries in a commutative ring $R$ is called weakly determinant liftable if there exists a matrix $B$ congruent to $A$ modulo $R\det(A)$ and $\det(B)=0$; if we can choose $B$ to be unimodular, then $A$ is called determinant liftable. If $A$ is extendable to an invertible $3\times 3$ matrix $A^+$, then $A$ is weakly determinant liftable. If $A$ is simply extendable (i.e., we can choose $A^+$ such that its $(3,3)$ entry is $0$), then $A$ is determinant liftable. We present necessary and/or sufficient criteria for $A$ to be (weakly) determinant liftable and we use them to show that if $R$ is a $\Pi_2$ ring in the sense of Part I (resp.\ is a pre-Schreier domain), then $A$ is simply extendable (resp.\ extendable) iff it is determinant liftable (resp.\ weakly determinant liftable). As an application we show that each $J_{2,1}$ domain (as defined by Lorenzini) is an elementary divisor domain.
\end{abstract}

\maketitle

\section{Introduction}\label{S1}

Let $R$ be a commutative ring with identity element $1$. For $n\in\mathbb{N}=\{1,2,\ldots\}$, let $\mathbb{M}_n(R)$ be the $R$-algebra of $n\times n$ matrices with entries in $R$. We say that $B,C\in\mathbb M_n(R)$ are congruent modulo an ideal $\mathfrak{i}$ of $R$ if all entries of $B-C$ belong to $\mathfrak{i}$, i.e., $B-C\in\mathbb M_n(\mathfrak{i})$. Let ${SL}_{n}(R):=\{I\in\mathbb M_n(R)|\det(I)=1\}$. For a free $R$-module $F$, let $Um(F)$ be the set of \textsl{unimodular} elements of $F$, i.e., of elements $v\in F$ for which there exists an $R$-linear map $L:F\rightarrow R$ such that $L(v)=1$. 

In this paper we study a unimodular matrix $A =\left[ 
\begin{array}{cc}
a & b \\ 
c & d\end{array}\right]\in Um\bigl(\mathbb M_2(R)\bigr)$. Recall that $A$ is called {\it extendable} if there exists $A^+=\left[ 
\begin{array}{ccc}
a & b & f \\ 
c & d & -e \\ 
-t & s & v\end{array}\right]\in SL_3(R)$ (see \cite{CPV1}, Def.\ 1.1); we call $A^+$ an {\it extension} of $A$. If we can choose $A^+$ such that $v=0$, then $A$ is called {\it simply extendable} and $A^+$ is called a {\it simple extension} of $A$.

\begin{definition}\label{def1}
We say that $A\in Um\bigl(\mathbb{M}_2(R)\bigr)$ is weakly determinant liftable if there exists $B\in \mathbb M_2(R)$ congruent to $A$ modulo $R\det(A)$ and $\det(B)=0$. If there exists such a matrix $B$ which is unimodular, then we remove the word `weakly', i.e., we say that $A$ is determinant liftable.
\end{definition}

If either $A$ is invertible or $\det(A)=0$, then $A$ is determinant liftable. Also, if $A^{\prime}\in\mathbb M_2(R)$ is equivalent to $A$, then $A^{\prime}$ is (weakly) determinant liftable iff $A$ is so. The following characterizations of determinant liftability are proved in Section \ref{S4}.

\begin{theorem}\label{TH1}
For $A=\left[ 
\begin{array}{cc}
a & b \\ 
c & d\end{array}\right]\in Um\bigl(\mathbb M_2(R)\bigr)$ the following statements are equivalent.

\medskip \textbf{(1)} The matrix $A$ is determinant liftable.

\smallskip \textbf{(2)} There exists $C\in Um\bigl(\mathbb{M}_{2}(R)\bigr) $ such that $A+\det(A)C\in Um\bigl(\mathbb{M}_{2}(R)\bigr)$ and $\det(C)=\det\bigl(A+\det(A)C\bigr)=0$. 

\smallskip \textbf{(3)} There exists $(x,y,z,w)\in R^4$ such that $ax+by+cz+dw=1$ and $xw-yz=0$.

\smallskip \textbf{(4)} There exists $C\in \mathbb{M}_2(R)$ such that $\det(C)=\det\bigl(A+\det(A)C\bigr)=0$.
\end{theorem}

Definition \ref{def1} is motivated by the following implications proved in Section \ref{S5}.

\begin{theorem}\label{TH2} For $A\in Um\bigl(\mathbb M_2(R)\bigr)$ the following properties hold.

\medskip \textbf{(1)} If $A$ is simply extendable, then $A$ is determinant liftable.

\smallskip \textbf{(2)} If $A$ is extendable, then $A$ is weakly determinant liftable.
\end{theorem}

The converses of Theorem \ref{TH2} do not hold in general (see Example \ref{EX3.5}).

Recall from \cite{CPV1}, Def.\ 1.2(1) that $R$ is called a {\it $\Pi_2$ ring} if each matrix in $Um\bigl(\mathbb M_2(R)\bigr)$ of zero determinant is extendable, equivalently it is simply extendable by \cite{CPV1}, Lem.\ 4.1(1). We will use stable ranges and pre-Schreier domains as recalled in \cite{CPV1}, Def.\ 1.5 and Sect.\ 2. Each pre-Schreier domain is a $\Pi_2$ ring (see \cite{CPV1}, paragraph after Thm.\ 1.4). Recall from \cite{CPV1}, Def.\ 1.2(3) that $R$ is called an {\it $SE_2$ ring} if each matrix in $Um\bigl(\mathbb M_2(R)\bigr)$ is simply extendable. In particular, an $SE_2$ ring is a $\Pi_2$ ring.

The next two theorems are proved in Section \ref{S7}.

\begin{theorem}\label{TH3}
The ring $R$ is a $\Pi_2$ ring iff the simply extendable and determinant liftable properties on a matrix in $Um\bigl(\mathbb M_2(R)\bigr)$ are equivalent.
\end{theorem}

From Theorem \ref{TH3} and \cite{CPV1}, Thm.\ 1.6 we get directly the following result.

\begin{corollary}\label{C1}
If $R$ is a $\Pi_2$ ring with $sr(R)\le 2$, then the simply extendable, extendable and determinant liftable properties on a matrix in $Um\bigl(\mathbb M_2(R)\bigr)$ are equivalent.
\end{corollary}

As $SE_2$ rings are $\Pi_2$ rings, from Theorem \ref{TH3} we get directly the following result.

\begin{corollary}\label{C2}
The ring $R$ is an $SE_2$ ring iff it is a $\Pi_2$ ring with the property that each matrix in $Um\bigl(\mathbb M_2(R)\bigr)$ with nonzero determinant is determinant liftable.
\end{corollary}

A matrix $N\in\mathbb M_2(R)$ is called {\it non-full} if it is a product $\left[ 
\begin{array}{c}
l \\ 
m\end{array}\right] \left[ 
\begin{array}{cc}
o & q\end{array}\right] $ with $(l,m,o,q)\in R^4$; so $\det(N)=0$. Recall that an integral domain $R$ is a pre-Schreier domain iff each matrix in $\mathbb{M}_{2}(R)$ of zero determinant is non-full by \cite{MR}, Lem.\ 1.

\begin{theorem}\label{TH4}
Assume $R$ is such that each zero determinant matrix in $\mathbb M_2(R)$ is non-full (e.g., $R$ is a product of pre-Schreier domains). Then a unimodular matrix $A\in Um\bigl(\mathbb{M}_{2}(R)\bigr)$ is extendable iff it is weakly determinant liftable.
\end{theorem}

Note that $R$ of Theorem \ref{TH4} is a $\Pi_2$ ring by \cite{CPV1}, Thm.\ 1.4. Based on this, from Theorems \ref{TH2} and \ref{TH4} and Corollary \ref{C1} we get directly the following result.

\begin{corollary}\label{C3}
Assume $R$ is such that $sr(R)\le 2$ and each zero determinant matrix in $\mathbb M_2(R)$ is non-full (e.g., $R$ is a product of pre-Schreier domains of stable range at most $2$). Then the simply extendable, extendable, determinant liftable and weakly determinant liftable properties on a matrix in $Um\bigl(\mathbb M_2(R)\bigr)$ are equivalent.
\end{corollary}

\begin{example}\label{EX1}
\normalfont 
Let $R$ be a pre-Schreier domain such that $sr(R)=3$ and there exists $A\in Um\bigl(\mathbb M_2(R)\bigr)$ which is extendable but not simply extendable (e.g., $R=K[X,Y]$ with $K$ a subfield of $\mathbb R$, see \cite{CPV1}, Ex.\ 6.1). Then $A$ is weakly determinant liftable by Theorem \ref{TH2}(2) but it is not determinant liftable by Theorem \ref{TH3}. So the inequalities in Corollaries \ref{C1} and \ref{C3} are optimal.
\end{example}

Recall that $R$ is called an {\it elementary divisor ring} if for each $(j,n)\in\mathbb N^2$ with $j\le n$, every matrix of size either $j\times n$ or $n\times j$ with entries in $R$ admits diagonal reduction, i.e., is equivalent to a matrix whose off diagonal entries are $0$ and whose diagonal entries $a_{1,1},\ldots,a_{j,j}$ are such that $a_{i,i}$ divides $a_{i+1,i+1}$ for all $i\in \{1,\ldots ,j-1\}$. Also, recall that $R$ is a {\it Hermite ring} in the sense of Kaplansky if $R^2=RUm(R^2)$, equivalently
if each $1\times 2$ matrix with entries in $R$ admits diagonal reduction. Clearly, each elementary divisor ring is a Hermite ring. Moreover, each elementary divisor ring is an $SE_2$ ring (as defined above) by \cite{CPV1}, Prop.\ 1.3.

Lorenzini introduced $3$ classes of rings that are `between' elementary divisor rings and Hermite rings (see \cite{lor}, Prop.\ 4.11). We define the first class, $J_{2,1}$, as follows.

\begin{definition}\label{def2}
We say that $R$ is:

\medskip \textbf{(1)} a ${WJ}_{2,1}$ ring if for each $(a,b,c,d)\in Um(R^4)$ and every $(\Psi,\Delta)\in R^2$, there exists $(x,y,z,w)\in R^{4}$ such that $ax+by+cx+dw=\Psi$ and $xw-yz=\Delta$;

\smallskip \textbf{(2)} a $J_{2,1}$ ring if it is a Hermite ring and a ${WJ}_{2,1}$ ring.
\end{definition}

The above definition of a $J_{2,1}$ ring is equivalent to the one in \cite{lor}, Def.\ 4.6 (see Proposition \ref{P1}). For $J_{n,r}$ rings with $(n,r)\in\mathbb N^2$ and $n>r$ see \cite{lor}, Def.\ 4.6. Lorenzini shows that an elementary divisor ring is a $J_{n,1}$ ring for every integer $n>1$ (see \cite{lor}, Prop.\ 4.8) and Fresnel shows that an Euclidean domain is a $J_{n,r}$ ring for all $(n,r)\in\mathbb N^2$ with $n>r$ (see \cite{fre}, Thm.\ 1.1). Theorem \ref{TH5} below, proved in Section \ref{S8} based on Theorem \ref{TH3}, solves a problem posed by Lorenzini (see \cite{lor}, p.\ 618) and Fresnel (see \cite{fre}, Subsect.\ 3.1) for the case of $\Pi_2$ rings.

\begin{theorem}\label{TH5}
Let $R$ be a ${WJ}_{2,1}$ ring. Then the following properties hold.

\medskip \textbf{(1)} Each matrix $A\in Um\bigl(\mathbb M_2(R)\bigr)$ is determinant liftable.

\smallskip \textbf{(2)} Assume $R$ is also a Hermite ring (i.e., $R$ is a $J_{2,1}$ ring). Then $R$ is an elementary divisor ring iff it is a $\Pi_2$ ring. 

\smallskip \textbf{(3)} Let $R$ be a $J_{2,1}$ domain. Then $R$ is an elementary divisor domain.
\end{theorem}

Therefore the constructions for an arbitrary commutative ring performed in \cite{lor}, Ex.\ 4.10 for $n=2$ or in Ex.\ 3.5 always produce rings which are either elementary divisor rings or are not $\Pi_2$ rings (in particular, the integral domains produced are elementary divisor domains). 

By combining Theorem \ref{TH5}(2) with \cite{lor}, Prop.\ 4.8 we get directly the following implication between these classes of rings.

\begin{corollary}\label{C4}
Let $R$ be a $J_{2,1}$ ring. If $R$ is a $\Pi_2$ ring, then $R$ is a $J_{n,1}$ ring for each integer $n>1$.
\end{corollary}

The commutative $R$-algebras associated to $A$ and required in the proofs of the above theorems are introduced in Section \ref{S2}; their smoothness properties are presented in Theorem \ref{TH6} of Section \ref{S3}. Section \ref{S6} proves criteria for weakly determinant liftability. Section \ref{S9} studies rings $R$ for which all $A\in Um\bigl(\mathbb M_2(R)\bigr)$ are (weakly) determinant liftable. Section \ref{S10} proves criteria for determinant liftability via completions. Section \ref{S11} uses Picard groups to refine Theorem \ref{TH6}(7) when $\det(A)=0$. Part III of this series of articles \cite{CPV2} will contain applications to Hermite rings.

\section{Five Algebras}\label{S2}

For a commutative $R$-algebra $S$, let $N(S)$, $Z(S)$, $U(S)$, $\textup{Pic\,}(S)$ be the nilradical, the set of
zero divisors, the multiplicative group of units and the Picard group (respectively) of $S$. Let $\textup{Spec\,} S$ be the spectrum of $S$. Let $\textup{Max\,} S$ be the set of maximal prime ideals of $S$. For $f\in S$, let $(f):=Sf$ and we often abuse the notation by denoting $\bar{f}:=f+\mathfrak{i}\in S/\mathfrak{i}$ for several ideals $\mathfrak{i}$ of $S$. 

To $A=\left[ 
\begin{array}{cc}
a & b \\ 
c & d\end{array}\right]\in Um\bigl(\mathbb M_2(R)\bigr)$ we attach five commutative $R$-algebras as follows. Let $X,Y,Z,W,V$ and $T$ be variables. The first four $R$-algebras are
\begin{equation*}
\begin{array}{c}
\mathcal U=\mathcal U_A:=R[X,Y,Z,W]/(1-aX-bY-cZ-dW),\;\;\;\;\;\;\;\;\;\;\;\;\;\;\;\;\;\;\;\;\;\;\;\;\;\;\;\;\;\;\;\;\;\;\;\;\;\\
\;\;\,\mathcal E=\mathcal E_A:=R[X,Y,Z,W,V]/\bigl(1-aXW-bXZ-cYW-dYZ-(ad-bc)V\bigr),\\
\mathcal X=\mathcal X_A:=R[X,Y,Z,W]/(1-aXW-bXZ-cYW-dYZ),\,\;\;\;\;\;\;\;\;\;\;\;\;\;\;\;\;\;\;\;\;\;\;\;\;\\
\mathcal D=\mathcal D_A:=R[X,Y,Z,W]/(1-aX-bY-cZ-dW,XW-YZ).\!\;\;\;\;\;\;\;\;\;\;\;\;\;\;\;\;\;\;\;\;\end{array}
\end{equation*}
The polynomial
$$\Phi=\Phi_A(X,Y,Z,W):=1-aX-bY-cZ-dW+(ad-bc)(XW-YZ)\in R[X,Y,Z,W]$$
admits several natural decompositions of the form $e_{1,1}e_{2,2}-e_{1,2}e_{2,1}$, such as $\break\Phi=(1-aX-cZ)(1-bY-dW)-(aY+cW)(bX+dZ)$. 
The fifth $R$-algebra is
\begin{equation*}
\mathcal W=\mathcal W_A:=R[X,Y,Z,W]/(\Phi).
\end{equation*}

If $c=0$, then 
\begin{equation}\label{EQ2}
\begin{array}{c}
\mathcal W=R[X,Y,Z,W]/\bigl((1-aX)(1-dW)-Y(b+adZ)\bigr)\\ 
\;\;\;\;\;\;\;\;\;\;\;\;\;\;\;\;\;\;=R[X,Y,Z,W]/\bigl((1-aX)(1-bY-dW)-aY(bX+dZ)\bigr).\end{array}
\end{equation}


We define an arrow diagram of $R$-algebra homomorphisms 

\begin{equation}\label{EQ3}
\xymatrix@R=10pt@C=21pt@L=2pt{
\mathcal W \ar[r] & \mathcal D\ar[d]^{\rho} & \mathcal U\ar[l]\\
\mathcal E \ar[r] & \mathcal X\\}
\end{equation}
as follows. The horizontal homomorphisms are surjections defined by $R$-algebra isomorphisms
$\mathcal W/(\bar{X}\bar{W}-\bar{Y}\bar{Z})\cong\mathcal D\cong\mathcal U/(\bar{X}\bar{W}-\bar{Y}\bar{Z})$
and $\mathcal X\cong\mathcal E/(\bar{V})$ given by the Third Isomorphism Theorem, and $\rho $ is defined by mapping
$\bar{X}$, $\bar{Y}$, $\bar{Z}$, $\bar{W}$ to $\bar{X}\bar{W}$, $\bar{X}\bar{Z}$, $\bar{Y}\bar{W}$, $\bar{Y}\bar{Z}$ (respectively). Equation (\ref{EQ2}) and Diagram (\ref{EQ3}) encode many applications in this Part II and the sequel Part III \cite{CPV2}. 

The $R$-algebra homomorphism $R[X,Y,Z,W]\rightarrow \mathcal D[T,T^{-1}]\otimes_{\mathcal D,\rho} \mathcal X$ that maps $X$, $Y$, $Z$, and $W$ to $T\otimes {\bar X}$, $T\otimes {\bar Y}$, $T^{-1}\otimes {\bar Z}$, and $T^{-1}\otimes {\bar W}$ (respectively) maps $1-aXW-bXZ-cYW-dYZ$ to $1\otimes (1-a{\bar X}{\bar W}-b{\bar X}{\bar Z}-c{\bar Y}{\bar W}-d{\bar Y}{\bar Z})=1\otimes 0=0$ and hence it induces an $R$-algebra homomorphism
$$\chi:\mathcal X\rightarrow\mathcal D[T,T^{-1}]\otimes_{\mathcal D,\rho} \mathcal X$$
that maps $\bar{X}$, $\bar{Y}$, $\bar{Z}$, and $\bar{W}$ to $T\otimes \bar{X}$, $T\otimes \bar{Y}$, $T^{-1}\otimes \bar{Z}$, and $T^{-1}\otimes\bar{W}$ (respectively). The composite $R$-algebra homomorphism $\chi\circ\rho:\mathcal D\rightarrow\mathcal D[T,T^{-1}]\otimes_{\mathcal D,\rho} \mathcal X$ maps $\bar{X}$, $\bar{Y}$, $\bar{Z}$, and $\bar{W}$ to $1\otimes \bar{X}\bar{W}={\bar X}\otimes 1$, $1\otimes \bar{X}\bar{Z}={\bar Y}\otimes 1$, $1\otimes \bar{Y}\bar{W}={\bar Z}\otimes 1$, $1\otimes \bar{Y}\bar{Z}={\bar W}\otimes 1$ (respectively) and hence $\chi$ is a $\mathcal D$-algebra homomorphism.

For two $R$-algebras $S$ and $S^{\prime}$, we consider the set
$$\textup{Hom}_{R}(S,S^{\prime}):=\{\varrho:S\rightarrow S^{\prime}|\varrho\;\textup{is an}\; R\textup{-algebra homomorphism}\}.$$ 
The next two lemmas form the first justification for introducing these $R$-algebras.

\begin{lemma}\label{L1}
Let $A=\left[ 
\begin{array}{cc}
a & b \\ 
c & d\end{array}\right]\in Um\bigl(\mathbb M_2(R)\bigr)$ and the $R$-algebras $\mathcal U$, $\mathcal E$, $\mathcal X$, $\mathcal D$ and $\mathcal W$ be as above. Then the following properties hold.

\medskip
{\bf (1)} For a commutative $R$-algebra $S$ we have functorial bijections of sets
\begin{equation}\label{EQ3.a}
\textup{Hom}_{R}(\mathcal U,S)\cong\{(x,y,z,w)\in S^4|ax+by+cz+dw=1\},\;\;\;\;\;\;\;\;\;\;\;\;\;\;\;\;\;\;\;\;\;\;\;\;\;
\end{equation}
\begin{equation}\label{EQ3.b}
\;\;\;\textup{Hom}_{R}(\mathcal E,S)\cong\{(x,y,z,w,v)\in S^5|axw+bxz+cyw+dyz+(ad-bc)v=1\},
\end{equation}
\begin{equation}\label{EQ3.c}
\textup{Hom}_{R}(\mathcal X,S)\cong\{(x,y,z,w)\in S^4|axw+bxz+cyw+dyz=1\},\;\;\;\;\;\;\;\;\;\;\;\;\;\;\;
\end{equation}
\begin{equation}\label{EQ3.e}
\textup{Hom}_{R}(\mathcal D,S)\cong\{(x,y,z,w)\in S^4|ax+by+cz+dw=1\;\textup{and}\;xw=yz\},\;\;\;
\end{equation}
\begin{equation}\label{EQ3.d}
\textup{Hom}_{R}(\mathcal W,S)\cong\{(x,y,z,w)\in S^4|\Phi(x,y,z,w)=0\}.\;\;\;\;\;\;\;\;\;\;\;\;\;\;\;\;\;\;\;\;\;\;\;\;\;\;\;\;\;\;\;\;\;\;\;\;
\end{equation}

\smallskip
{\bf (2)} For a commutative $R$-algebra $S$ we have a functorial diagram of sets
$$\xymatrix@R=10pt@C=21pt@L=10pt{
\{(x,y,z,w)\in S^4|axw+bxz+cyw+dyz=1\} \ar[r]\ar[d]_{\varrho(S)} & \textup{Hom}_R(\mathcal X,S)\ar[d]^{\rho(S)}\\
\{(x,y,z,w)\in S^4|ax+by+cz+dw=1,\; xw=yz\} \ar[r] & \textup{Hom}_R(\mathcal D,S),\\}$$
where the horizontal arrows are identifications given by Equations (\ref{EQ3.c}) and (\ref{EQ3.e}), $\rho(S)$ maps $h\in\textup{Hom}_R(\mathcal X,S)$ to $h\circ\rho\in\textup{Hom}_R(\mathcal D,S)$, and $\varrho(S)$ is defined by the rule $(x,y,z,w)\mapsto (xw,xz,yw,yz)$.\footnote{Thus the $R$-algebra homomorphism $\rho:\mathcal D\rightarrow\mathcal X$, as a morphism of functors from the category of $R$-algebras to the category of sets, represents decompositions of the unimodular matrix $\left[ 
\begin{array}{cc}
\bar X & \bar Y\\ 
\bar W & \bar Z\end{array}\right]\in Um\bigl(\mathbb M_2(\mathcal D)\bigr)$ as a product $\left[ 
\begin{array}{c}
\bar X_{\mathcal X}\\ 
\bar Y_{\mathcal X}\end{array}\right]\left[ 
\begin{array}{cc}
\bar W_{\mathcal X} & \bar Z_{\mathcal X}\end{array}\right]$, the lower right index ${\mathcal X}$ emphasizing variables for $\mathcal X$.\label{foo1}} 

{\bf (3)} Let $S$ be a commutative $R$-algebra. Let $\varrho(S)$ be as in part (2). Let $(x^{\prime},y^{\prime},z^{\prime},w^{\prime})\in S^4$ be such that $ax^{\prime}+by^{\prime}+cz^{\prime}+dw^{\prime}=1$, $x^{\prime}w^{\prime}=y^{\prime}z^{\prime}$, and $\{x^{\prime},y^{\prime},z^{\prime},w^{\prime}\}\cap U(S)\neq\emptyset$. Then there exists a quadruple $\upsilon=(x,y,z,w)\in S^4$ such that $axw+bxz+cyw+dyz=1$ and $\varrho(S)(\upsilon)=(x^{\prime},y^{\prime},z^{\prime},w^{\prime})$. Moreover, for each such quadruple $\upsilon$, the function $\varrho(S)_{\upsilon}: U(S)\times \{\upsilon\}\rightarrow \varrho(S)^{-1}(x^{\prime},y^{\prime},z^{\prime},w^{\prime})$ defined by $\varrho(S)_{\upsilon}(u,\upsilon):=(ux,uy,u^{-1}z,u^{-1}w)$ is a bijection.
\end{lemma}

\begin{proof}
Let $\textup{H}_S:=\{\varrho\in\textup{Hom}_{R}(R[X,Y,Z,W],S)|1-aX-bY-cZ-dW\in\textup{Ker}(\varrho)\}$. Factorization of homomorphisms theorem gives a functorial bijection
$\textup{Hom}_{R}(\mathcal U,S)\cong\textup{H}_S$. Similarly, the universal property of polynomial rings gives a functorial bijection 
$\textup{Hom}_{R}(R[X,Y,Z,W],S)\cong S^4$ defined by the rule $\varrho\mapsto\bigl(\varrho(X),\varrho(Y),\varrho(Z),\varrho(W)\bigr)$; this induces a functorial bijection $\textup{H}_S\cong\{(x,y,z,w)\in S^4|ax+by+cz+dw=1\}$. The last two sentences imply that Equation (\ref{EQ3.a}) holds.
Equations (\ref{EQ3.b}) to (\ref{EQ3.e}) are proved similarly. So part (1) holds. 

Part (2) follows from the definition of $\rho:\mathcal D\rightarrow\mathcal X$. 

For part (3), to fix the ideas we can assume that $x^{\prime}\in U(S)$; so $w^{\prime}=(x^{\prime})^{-1}y^{\prime}z^{\prime}$. As $\upsilon$ we can take $\bigl(1,(x^{\prime})^{-1}z^{\prime},y^{\prime},x^{\prime})\bigr)$. If $\varrho(S)_{\upsilon}(x_1,y_1,z_1,w_1)=(x^{\prime},y^{\prime},z^{\prime},w^{\prime})$, then $x_1w_1=xw=x^{\prime}\in U(S)$, $x_1z_1=xz=y^{\prime}$, $y_1w_1=yw=z^{\prime}$ and $y_1z_1=yz=w^{\prime}$. For $u\in U(S)$, $\varrho(S)_{\upsilon}(u,\upsilon)$ is equal to $(x_1,y_1,z_1,w_1)$ iff $u:=ww_1^{-1}$. So $\varrho(S)_{\upsilon}$ is a bijection and part (3) holds.\end{proof}

Recall that a retraction of an $R$-algebra $S$ is an $R$-algebra homomorphism $S\rightarrow R$. 

\begin{lemma}\label{L2}
For $A=\left[ 
\begin{array}{cc}
a & b \\ 
c & d\end{array}\right]\in Um\bigl(\mathbb M_2(R)\bigr)$ the following statements are equivalent.

\medskip
{\bf (1)} The matrix $A$ is extendable (resp.\ simply extendable).

\smallskip
{\bf (2)} There exists $(x,y,z,w,v)\in R^5$ (resp.\ $(x,y,z,w)\in R^4$) such that we have $axw+bxz+cyw+dyz+(ad-bc)v=1$ (resp.\ $axw+bxz+cyw+dyz=1$).

\smallskip
{\bf (3)} The $R$-algebra $\mathcal E$ (resp.\ $\mathcal X$) has a retraction.
\end{lemma}

\begin{proof}
The equivalence $(3)\Leftrightarrow (2)$ follows from Lemma \ref{L1}, Equation (\ref{EQ3.b}) (resp.\ (\ref{EQ3.c})) applied to the $R$-algebra $R$. By very definition, the matrix $A$ is extendable (resp.\ simply extendable) iff there exists $(x,y,z,w,v)\in R^5$ (resp.\ $(x,y,z,w)\in R^4$) such that the matrix $\left[ 
\begin{array}{ccc}
a & b & y\\ 
c & d & -x\\
-z & w & v
\end{array}\right]$ (resp.\ $\left[ 
\begin{array}{ccc}
a & b & y\\ 
c & d & -x\\
-z & w & 0
\end{array}\right]$) has determinant $1$, i.e., the identity $axw+bxz+cywdyz+(ad-bc)v=1$ (resp.\ $axw+bxz+cyw+dyz=1$) holds; hence $(1)\Leftrightarrow (2)$. Thus the lemma holds.\end{proof}

\section{Smoothness properties}\label{S3}

The analogs of Lemma \ref{L2} for determinant liftability (see Theorem \ref{TH1} and  Lemma \ref{L1}, Equation (\ref{EQ3.e})) and for weakly determinant liftability (see Theorem \ref{TH9}(2) to (4) below and  Lemma \ref{L1}, Equation (\ref{EQ3.d})) require substantial extra work and even some algebraic geometry reviews which are carried out in this section.

For an $R$-algebra homomorphism $h:S\rightarrow S^{\prime}$ between commutative $R$-algebras and $f\in R$, let $h_f:S_f\rightarrow S^{\prime}_f$ be the localization of $h$ with respect to the multiplicative set $\mathcal M_f:=\{f^i|i\in\mathbb N\cup\{0\}\}$ (with $f^0:=1$ by convention). In particular, $R_f:=\mathcal M_f^{-1}R$ and $S_f\cong R_f\otimes_R S$ as $R_f$-algebras. Similarly, for an $S$-linear map $\lambda:M\rightarrow M'$ between $S$-modules, let $\lambda_f: M_f\rightarrow M'_f$ be its localization with respect to $\mathcal M_f$; we have a functorial $S_f$-linear isomorphism $M_f\cong R_f\otimes_R M$, where $R_f\otimes_R M$ is viewed as an $S_f$-module via the $R_f$-algebra isomorphism $S_f\cong R_f\otimes_R S$. If $f\in N(R)$, then $M_f$ is the zero module over the zero ring $S_f$.

Given $n\in\mathbb N\cup\{0\}$, for the definition of a smooth morphism $\textup{Spec\,} S\rightarrow \textup{Spec\,} R$ (equivalently, a smooth homomorphism $R\rightarrow S$) of relative dimension $n$ we refer to \cite{BLR}, Ch.\ 2, Sect.\ 2.2, Def.\ 3. The main example of a smooth homomorphism of relative dimension $n$ is the polynomial $R$-algebra homomorphism $R\rightarrow R[X_1,\ldots,X_n]$. Smoothness of relative dimension $n$ is a local notion in the Zariski topologies of either $\textup{Spec\,} S$ or $\textup{Spec\,} R$ and hence (i) the homomorphism $R\rightarrow S$ is smooth of relative dimension $n$ if $S$ is a symmetric $R$-algebra of a projective $R$-module of rank $n$ and (ii) the localization homomorphism $R\rightarrow R_f$ is smooth of relative dimension $0$. If $r\in\mathbb N\cup\{0\}$ and the homomorphisms $R\rightarrow S$ and $S\rightarrow S^{\prime}$ are smooth of relative dimensions $n$ and $r$ (respectively), then the composite homomorphism $R\rightarrow S^{\prime}$ is smooth of relative dimension $n+r$.

Recall that $\mathbb{G}_{m,R}$ is the commutative affine group scheme over $\textup{Spec\,} R$ defined by the commutative Hopf $R$-algebra $R[T,T^{-1}]$ whose $R$-algebra comultiplication homomorphism $\textup{co}_R:R[T,T^{-1}]\rightarrow R[T,T^{-1}]\otimes_R R[T,T^{-1}]$ is defined by $\textup{co}_R(T)=T\otimes T$ (the $R$-algebra homomorphisms that are the counit $R[T,T^{-1}]\rightarrow R$ and the coinverse $R[T,T^{-1}]\rightarrow R[T,T^{-1}]$ are uniquely determined by the comultiplication and they map $T$ to $1$ and respectively $T^{-1}$); so for each commutative $R$-algebra $S$, we have a functorial isomorphism of abstract commutative groups
$$\{\hbar:\textup{Spec\,} S\rightarrow \mathbb{G}_{m,R}|\hbar\;\textup{is a morphism of schemes over}\;\textup{Spec\,} R\}\cong U(S).$$

Recall that a morphism $\textup{Spec\,} S\rightarrow\textup{Spec\,} R$ with a (left) action of $\mathbb G_{m,R}$ on it defined by a morphism $\mathcal A:\mathbb G_{m,R}\times_{\textup{Spec\,} R} \textup{Spec\,} S\rightarrow \textup{Spec\,} S$ is a {\it $\textup{Spec\,} R$-torsor} (in the Zariski topology) {\it under} $\mathbb{G}_{m,R}$ if there exists $n\in\mathbb N$ and $(f_1,\ldots,f_n)\in Um(R^n)$ (so $\textup{Spec\,} R=\cup_{i=1}^n \textup{Spec\,} R_{f_i}$) such that for each $i\in\{1,\ldots,n\}$ there exists a $\textup{Spec\,} R_{f_i}$-isomorphism $\textup{Spec\,} S_{f_i}\rightarrow \mathbb G_{m,R_{f_i}}$ under which the action $\mathbb G_{m,R_{f_i}}\times_{\textup{Spec\,} R_{f_i}} \textup{Spec\,} S_{f_i}\rightarrow \textup{Spec\,} S_{f_i}$ becomes isomorphic to the product morphism $\mathbb G_{m,R_{f_i}}\times_{\textup{Spec\,} R_{f_i}} \mathbb G_{m,R_{f_i}}\rightarrow \mathbb G_{m,R_{f_i}}$ defined by the comultiplication $\textup{co}_{R_{f_i}}$. 

If we can take $n=1$, then $\textup{Spec\,} S\rightarrow\textup{Spec\,} R$ with the action of $\mathbb G_{m,R}$ on it is called a {\it trivial $\textup{Spec\,} R$-torsor under $\mathbb{G}_{m,R}$}. If $n=1$, then $f_1\in U(R)$ and it follows that there exist $R$-algebra isomorphisms $S\cong S_{f_1}\cong R_{f_1}[T,T^{-1}]\cong R[T,T^{-1}]$ and in particular that the $R$-algebra $S$ has a retraction. Conversely, if the $R$-algebra $S$ has a retraction $h:S\rightarrow R$, then the natural composite morphism 
$$\mathbb{G}_{m,R}\cong \mathbb{G}_{m,R}\times_{\textup{Spec\,} R} \textup{Spec\,} R\rightarrow \mathbb{G}_{m,R}\times_{\textup{Spec\,} R} \textup{Spec\,} S\rightarrow\textup{Spec\,} S$$ defined by restricting $\mathcal A$ to the closed subscheme $\textup{Spec\,} h:\textup{Spec\,} R\rightarrow \textup{Spec\,} S$ of $\textup{Spec\,} S$ is an isomorphism between $\textup{Spec\,} R$-schemes and its inverse is a $\textup{Spec\,} R$-isomorphism $\textup{Spec\,} S\rightarrow \mathbb G_{m,R}$ under which $\mathcal A$ becomes isomorphic to the product morphism $\mathbb G_{m,R}\times_{\textup{Spec\,} R} \mathbb G_{m,R}\rightarrow \mathbb G_{m,R}$ and hence we can take $n=1$ with $f_1:=1$.

If $S$ is a commutative $R$-algebra, the pullback of a $\textup{Spec\,} R$-torsor under $\mathbb{G}_{m,R}$ via the morphism $\textup{Spec\,} S\rightarrow\textup{Spec\,} R$ is a $\textup{Spec\,} S$-torsor under $\mathbb{G}_{m,S}$. 

For more details related to the last two paragraphs see \cite{BLR}, Ch.\ 6, Sect.\ 5.4 or \cite{M}, Ch.\ III, Sect.\ 4 which use the $\textup{fppq}$ and the flat (respectively) topology; see \cite{M}, Ch.\ III, Prop.\ 4.9 for why we can restrict to the Zariski topology. 

As $A=\left[ 
\begin{array}{cc}
a & b \\ 
c & d\end{array}\right]\in Um\bigl(\mathbb M_2(R)\bigr)$, the $R$-linear map $\Lambda_A:R^4\rightarrow R$ defined by the rule $\Lambda_A(x,y,z,w)=ax+by+cz+dw$ is surjective. Thus 
\begin{equation}\label{EQP}
P=P_A:=\textup{Ker}(\Lambda_A)=\{(x,y,z,w)\in R^4|ax+by+cz+dw=0\}
\end{equation}
and its dual $P^*:=\{f:P\rightarrow R|f\;\textup{is}\;R\textup{-linear}\}$ are projective $R$-modules of rank $3$ with $P\oplus R\cong P^*\oplus R\cong R^4$. Also, $P\cong P^*$ by \cite{lam}, Ch.\ III, Sect.\ 6, Thm.\ 6.7 (1).

\begin{theorem}\label{TH6} 
For $A=\left[ 
\begin{array}{cc}
a & b \\ 
c & d\end{array}\right]\in Um\bigl(\mathbb M_2(R)\bigr)$, let $P$ be as above. The following properties hold for the $R$-algebras $\mathcal U$, $\mathcal E$, $\mathcal X$, $\mathcal D$, and $\mathcal W$ of Section \ref{S2}.

\medskip \textbf{(1)} For $f\in \{a,b,c,d\}$, the $R_{f}$-algebra $\mathcal U_f$ is isomorphic to a polynomial $R_f$-algebra in $3$ variables. In particular, the $R$-algebra $\mathcal U$ is smooth of relative dimension $3$.

\smallskip \textbf{(2)} The $R$-algebra $\mathcal U$ is isomorphic to the symmetric $R$-algebra of $P$.

\smallskip \textbf{(3)} The $R$-algebra $\mathcal D$ is smooth of relative dimension $2$.

\smallskip \textbf{(4)} The $R$-algebra $\mathcal E$ is smooth of relative dimension $4$.

\smallskip \textbf{(5)} The localization $\mathcal W_{\bigl(1-(ad-bc)({\bar X}{\bar W}-{\bar Y}{\bar Z})\bigr)}$ of $\mathcal W$ is a smooth
$R$-algebra of relative dimension $3$. In particular, if $ad-bc\in N(R)$, then the $R$-algebra $\mathcal W$ is smooth of relative dimension $3$.

\smallskip \textbf{(6)} We consider the morphism of $\textup{Spec\,} R$-schemes $\textup{Spec\,} \rho:\textup{Spec\,}\mathcal X\rightarrow \textup{Spec\,} \mathcal D$ defined by the $R$-algebra homomorphism $\rho:\mathcal D\rightarrow\mathcal X$ of Section \ref{S2} and the morphism
$$\textup{Spec\,} \chi:\mathbb{G}_{m,\mathcal D}\times _{\textup{Spec\,}\mathcal D}\textup{Spec\,}\mathcal X\rightarrow \textup{Spec\,}\mathcal X$$ 
of $\textup{Spec\,} D$-schemes defined by the $\mathcal D$-algebra homomorphism $\chi$ of Section \ref{S2}. Then this morphism of $\textup{Spec\,}\mathcal D$-schemes is a $\mathbb{G}_{m,\mathcal D}$-action on $\textup{Spec\,} \rho$.

\smallskip \textbf{(7)} Referring to part (6), $\textup{Spec\,} \rho$ with its $\mathbb{G}_{m,\mathcal D}$-action is a $\textup{Spec\,}\mathcal D$-torsor under $\mathbb{G}_{m,\mathcal D}$. In particular, the morphism $\textup{Spec\,} \rho$ is smooth of relative dimension $1$ and the $R$-algebra $\mathcal X$ is smooth of relative dimension $3$.

\smallskip \textbf{(8)} Assume $\det(A)=0$. Then for $f\in \{a,b,c,d\}$ 
the $R_{f}$-algebra $\mathcal D_{f}$ is a polynomial $R_f$-algebra in $2$ variables. Moreover, there exists a self-dual projective $R$-module $Q$ of rank $2$ such that the $R$-algebra $\mathcal D$ is isomorphic to the symmetric $R$-algebra of $Q$ and $Q\oplus R\cong P$ (thus $Q\oplus R^2\cong R^4$). 
\end{theorem}

\begin{proof}
(1) As $A \in Um\bigl(\mathbb M_2(R)\bigr)$, we have $\textup{Spec\,} R=\cup _{f\in \{a,b,c,d\}}\textup{Spec\,} R_{f}$. For $f=a$, as $a$ is invertible in $R_a$, the $R_a$-algebra $\mathcal U_a$ is easily seen to be isomorphic to $R_a[Y,Z,W]$. The other three cases $f\in\{b,c,d\}$ are similar.

(2) We consider the $R$-algebra
$$\mathcal U^{\prime}:=R[X^{\prime},Y^{\prime},Z^{\prime},W^{\prime}]/(aX^{\prime}+bY^{\prime}+cZ^{\prime}+dW^{\prime}).$$
Let $(a^{\prime},b^{\prime},c^{\prime},d^{\prime})\in R^4$ be such that $\Lambda_A(a^{\prime},b^{\prime},c^{\prime},d^{\prime})=1$.  The $R$-algebra isomorphism $R[X,Y,Z,W]\rightarrow R[X^{\prime},Y^{\prime},Z^{\prime},W^{\prime}]$ that maps $X$ to $a^{\prime}-X^{\prime}$, $Y$ to $b^{\prime}-Y^{\prime}$, $Z$ to $c^{\prime}-Z^{\prime}$, and $W$ to $d^{\prime}-W^{\prime}$ also maps $1-aX-bY-cZ-dW$ to $aX^{\prime}+bY^{\prime}+cZ^{\prime}+dW^{\prime}$ and hence it induces an $R$-algebra isomorphism $\mathcal U\cong\mathcal U^{\prime}$. 
We consider the $R$-module 
$$M_3:=(X^{\prime},Y^{\prime},Z^{\prime},W^{\prime})/\bigl(aX^{\prime}+bY^{\prime}+cZ^{\prime}+dW^{\prime}+(X^{\prime},Y^{\prime},Z^{\prime},W^{\prime})^2\bigr)$$ 
of quotient of ideals of $R[X^{\prime},Y^{\prime},Z^{\prime},W^{\prime}]$; it is isomorphic to the quotient $R$-module 
$$[RX^{\prime}\oplus RY^{\prime}\oplus RZ^{\prime}\oplus RW^{\prime}]/R(aX^{\prime},bY^{\prime},cZ^{\prime},dW^{\prime})$$ 
whose dual is naturally isomorphic to $P$. Thus the $R$-module $M_3$ is isomorphic to $P^*$ and hence also to $P$. We consider the $R$-linear map $l:M_3\rightarrow\mathcal U^{\prime}$ that maps $\star+\mathfrak I$ to $\star+\mathfrak J$ for each $\star\in\{X^{\prime},Y^{\prime},Z^{\prime},W^{\prime}\}$, where $\mathfrak I$ and $\mathfrak J$ are the ideals $\bigl(aX^{\prime}+bY^{\prime}+cZ^{\prime}+dW^{\prime}+(X^{\prime},Y^{\prime},Z^{\prime},W^{\prime})^2\bigr)$ and $(aX^{\prime}+bY^{\prime}+cZ^{\prime}+dW^{\prime})$ of $R[X^{\prime},Y^{\prime},Z^{\prime},W^{\prime}]$ (respectively). If $\mathfrak S$ is the symmetric algebra of the $R$-module $M_3$, then there exists a unique $R$-algebra homomorphism $\mathfrak l:\mathfrak S\rightarrow\mathcal U^{\prime}$ that extends $l$. To check that $\mathfrak l$ is an isomorphism it suffices to show that the $R_f$-algebra homomorphism $\mathfrak l_f:\mathfrak S_f\rightarrow\mathcal U_f^{\prime}$ is an isomorphism for each $f\in\{a,b,c,d\}$. For $f=a$, as $a$ is invertible in $R_a$, the $R_a$-algebra homomorphism $\mathfrak l_a$ is naturally identified with the identity automorphism of $R_a[Y,Z,W]$. The other three cases $f\in\{b,c,d\}$ are similar. We conclude that there exists an $R$-algebra isomorphism $\mathcal U\rightarrow\mathfrak S$ with $M_3\cong P$, thus part (2) holds.

(3) We consider the matrix 
\begin{equation*}
N_A:=\left[ 
\begin{array}{cccc}
a & b & c & d \\ 
W & -Z & -Y & X\end{array}\right]
\in\mathbb M_2(R[X,Y,Z,W]).\end{equation*} For $\mathfrak{n}\in\textup{Max\,} (R[X,Y,Z,W])$ such that $\{XW-YZ,1-aX-bY-cZ-dW\}\subset\mathfrak n$. Let $\textup{r}_{\mathfrak n}\in\{0,1,2\}$ be the rank of $N_A$ modulo $\mathfrak n$ and let $\kappa
:=R[X,Y,Z,W]/\mathfrak{n}$. We show that the assumption that $\textup{r}_{\mathfrak n}\le 1$ leads to a contradiction. As $N_A$ modulo $\mathfrak{n}$ has unimodular rows, we have $\textup{r}_{\mathfrak n}=1$ and thus there exists $\alpha \in R[X,Y,Z,W]\setminus\mathfrak{n}$ such that $(a,b,c,d)+\alpha (W,-Z,-Y,X)\in \mathfrak{n}^{4}$. So $2\alpha (XW-YZ)$
is congruent to $(aX+bY+cZ+dW)$ modulo $\mathfrak{n}$, and thus, as $\{XW-YZ,1-aX-bY-cZ-dW\}\subset\mathfrak{n}$, it is congruent to both $0$ and $-1$ modulo $\mathfrak{n}$, a contradiction. Therefore $\textup{r}_{\mathfrak n}=2$.

We denote the differential forms operator by $\delta$ in order to avoid confusion with the element $d\in R$. As $\textup{r}_{\mathfrak n}=2$, the two differential forms $\delta(XW-YZ)=W\delta X-Z\delta Y-Y\delta Z+X\delta W$ and $\delta(1-aX-bY-cZ-dW)=-a\delta X-b\delta Y-c\delta Z-d\delta W$ of the free $R[X,Y,Z,W]$-module $\Omega^1_{R[X,Y,Z,W]/R}=\oplus_{\star\in\{X,Y,Z,W\}} R[X,Y,Z,W]\delta\star$ of relative differential forms (of degree $1$) of $R[X,Y,Z,W]$ over $R$ (see \cite{BLR}, Ch.\ 2, Sect.\ 2.1, p. 31) are linearly independent in $\Omega^1_{R[X,Y,Z,W]/R}\otimes_{R[X,Y,Z,W]} \kappa$. From this and Jacobian Criterion (see the implication $(a)\Leftrightarrow (d)$ in \cite{BLR}, Ch.\ 2, Sect.\ 2.2, Prop.\ 7) it follows that the morphism $\textup{Spec\,}\mathcal D\rightarrow\textup{Spec\,} R$ is smooth of relative dimension $4-2=2$ at $\mathfrak n/(XW-YZ,1-aX-bY-cZ-dW)\in\textup{Spec\,}\mathcal D$. From this and \cite{BLR}, Ch.\ 2, Sect.\ 2.2, Prop.\ 11 it follows that $\textup{Spec\,}\mathcal D\rightarrow\textup{Spec\,} R$ is smooth of relative dimension $2$ at all points of $\textup{Spec\,} D$ contained in an open subset of $\textup{Spec\,}\mathcal D$ that contains $\textup{Max\,}\mathcal D$ and thus at all points of $\textup{Spec\,}\mathcal D$. Thus part (3) holds by the very definitions (see \cite{BLR}, Ch.\ 2, Sect.\ 2.2, Def.\ 3).

(4) Similar to the last two paragraphs, it suffices to 
show that for 
$$\Theta(X,Y,Z,W,V):=1-aXW-bXZ-cYW-dYZ-(ad-bc)V\in R[X,Y,Z,W,V],$$ 
and its five partial derivatives $\Theta_X$, $\Theta_Y$, $\Theta_Z$, $\Theta_W$, and $\Theta_V$, the differential form 
$\delta\Theta=\Theta_X\delta X+\Theta_Y\delta Y+\Theta_Z\delta_Z+\Theta_W\delta W+\Theta_V\delta V\in\Omega^1_{R[X,Y,Z,W,V]/R}$ is nonzero modulo each maximal ideal of $R[X,Y,Z,W,V]$ that contains $\Theta$. So it suffices to show that $\Theta$ and its five partial derivatives generate $R[X,Y,Z,W,V]$, which follows from the identity $1=\Theta-V\Theta_V-X\Theta_X-Y\Theta_Y$. So part (4) holds.

(5) Similarly to the last paragraph, part (5) follows from the fact that we have an identity $1-(ad-bc)(XW-YZ)=\Phi-X\Phi_X-Y\Phi_Y-Z\Phi_Z-W\Phi_W$. Note that if $ad-bc\in N(R)$, then $1-(ad-bc)({\bar X}{\bar W}-{\bar Y}{\bar Z})\in U(\mathcal W)$ and thus the $R$-algebra localization homomorphism $\mathcal W\rightarrow\mathcal W_{\bigl(1-(ad-bc)({\bar X}{\bar W}-{\bar Y}{\bar Z})\bigr)}$ is an isomorphism.\footnote{The morphism $\textup{Spec\,}\mathcal W\rightarrow\textup{Spec\,} R$ is not smooth of relative dimension $3$ at precisely the points of the closed subscheme of $\textup{Spec\,} \mathcal W$ which is the closed subscheme of $\textup{Spec\,}\mathcal W_{ad-bc}$ defined by the zero locus ${\bar X}-(ad-bc)^{-1}d={\bar Y}-(ad-bc)^{-1}c={\bar Z}-(ad-bc)^{-1}b={\bar W}-(ad-bc)^{-1}a=0$. In particular, if $A$ is invertible, then this closed subscheme is isomorphic to $\textup{Spec\,}\bigl(R/N(R)\bigr)$.}

(6) It is easy to see that we have an identity of $\mathcal D$-algebra homomorphisms
$$(\textup{co}_{\mathcal D}\otimes 1_{\mathcal X})\circ\chi=(1_{\mathcal D[T,T^{-1}]}\otimes \chi) \circ\chi:\mathcal X\rightarrow \mathcal D[T,T^{-1}]\otimes_{\mathcal D} \mathcal D[T,T^{-1}]\otimes_{\mathcal D,\rho} \mathcal X,$$
with $1_{\bigstar}$ denoting the identity automorphism of the $R$-algebra $\bigstar$ and the comultiplication $\textup{co}_{\mathcal D}:\mathcal D[T,T^{-1}]\rightarrow \mathcal D[T,T^{-1}]\otimes_{\mathcal D} \mathcal D[T,T^{-1}]$ as in the beginning of this section. For instance, $(\textup{co}_{\mathcal D}\otimes 1_{\mathcal X})\circ\chi({\bar Z})=(\textup{co}_{\mathcal D}\otimes 1_{\mathcal X})(T^{-1}\otimes {\bar Z})=T^{-1}\otimes T^{-1}\otimes{\bar Z}$ is equal to $(1_{\mathcal D[T,T^{-1}]}\otimes \chi) \circ\chi({\bar Z})=(1_{\mathcal D[T,T^{-1}]}\otimes \chi)(T^{-1}\otimes{\bar Z})=T^{-1}\otimes T^{-1}\otimes{\bar Z}$. 

Also, $\chi:\mathcal X\rightarrow \mathcal D[T,T^{-1}]\otimes_{\mathcal D,\rho} \mathcal X$ has a retraction defined by the $\mathcal X$-algebra homomorphism $\mathcal D[T,T^{-1}]\otimes_{\mathcal D,\rho} \mathcal X\rightarrow\mathcal X$ that maps $T\otimes 1$ to $1$.

Part (6) follows from the last two paragraphs.



(7) As $1=a\bar X+b\bar Y+c\bar Z+d\bar W$ in $\mathcal D$, $\textup{Spec\,}\mathcal D=\cup_{\bar f\in\{\bar X,\bar Y,\bar Z,\bar W\}} \textup{Spec\,}\mathcal D_{\bar f}$. So it suffices to show that for $\bar f\in\{\bar X,\bar Y,\bar Z,\bar W\}$, the pullback $\mathcal D_{\bar f}$-torsor under $\mathbb G_{m,\mathcal D_{\bar f}}$, i.e., the morphism $\textup{Spec\,} \mathcal X_{\bar f}\rightarrow \textup{Spec\,}\mathcal D_{\bar f}$ defined by $\rho_{\bar f}$ with its induced action, is trivial. We only show this for $\bar f=\bar X$ as the other three cases are similar; based on the diagram of Lemma \ref{L1}(2), this case follows from the last part of Lemma \ref{L1}(3) applied to $S=\mathcal D_{\bar X}$.

(8) For the polynomial $R$-algebra part we only consider the case $f=a$ as the other three cases are similar. By eliminating $\bar{X}=a^{-1}(1-b\bar{Y}-c\bar{Z}-d\bar{W})\in\mathcal D_a$, we obtain an isomorphism 
\begin{equation*}
\mathcal D_{a}\cong R_{a}[Y,Z,W]/\bigl(YZ-a^{-1}W(1-bY-cZ-dW)\bigr),
\end{equation*}which via the change of variables $(Y_1,Z_1,W_1):=(aY+cW,Z+a^{-1}bW,W)$ over $R_a$ is isomorphic, as $ad=bc$, to $R_a[Y_1,Z_1,W_1]/a^{-2}(Y_1Z_1-W_1)\cong R_a[Y_1,Z_1]$ and hence to the symmetric $R_a$-algebra of $R_a^2$. From \cite{BCW}, Thm.\ (4.4) it follows that $\mathcal D$ is isomorphic to the symmetric $R$-algebra of a projective $R$-module $Q$. From \cite{BCW}, Lem.\ (4.6) applied to $R$ and $R_f$ with $f\in\{a,b,c,d\}$ it follows that $Q$ is uniquely determined up to isomorphism and that $Q_f\cong R_f^2$ for each $f\in\{a,b,c,d\}$. Thus $Q$ has rank $2$. So there exists an $R$-algebra homomorphism $\mathcal D\rightarrow R$; let $\zeta:=(a^{\prime},b^{\prime},c^{\prime},d^{\prime})\in R^4$ be such that $\Lambda_A(\zeta)=1$ and $a^{\prime}d^{\prime}=b^{\prime}c^{\prime}$ by  Lemma \ref{L1}, Equation (\ref{EQ3.e}). The $R$-algebra $\mathcal D$ is isomorphic to
$$R[X^{\prime},Y^{\prime},Z^{\prime},W^{\prime}]/(aX^{\prime}+bY^{\prime}+cZ^{\prime}+dW^{\prime},d^{\prime}X^{\prime}-c^{\prime}Y^{\prime}-b^{\prime}Z^{\prime}+a^{\prime}W^{\prime}-X^{\prime}W^{\prime}+Y^{\prime}Z^{\prime})$$
via the change of variables $(X^{\prime},Y^{\prime},Z^{\prime},W^{\prime}):=(a^{\prime},b^{\prime},c^{\prime},d^{\prime})+(X,Y,Z,W)$.

We consider the $R$-submodule $M_2:=R(a,b,c,d)+ R(d^{\prime},-c^{\prime},-b^{\prime},a^{\prime})$ of $R^4$. To study $M_2$ we introduce the nondegenerate bilinear form $\langle\,,\rangle:R^4\times R^4\rightarrow R$ on the $R$-module $R^4$ defined by the rule 
$$\langle(x_1,x_2,x_3,x_4), (y_1,y_2,y_3,y_4)\rangle:=\sum_{i=1}^4 x_iy_i.$$ We have identities $\langle \zeta,(a,b,c,d)\rangle=1$ and $\langle\zeta,(d^{\prime},-c^{\prime},-b^{\prime},a^{\prime})\rangle=0$. 
Let $\mathfrak m\in\textup{Max\,} R$; so $R/\mathfrak m$ is a field and $V_4:=(R/\mathfrak m)^4$ is a vector space over it of dimension $4$. The image $V_2$ of $M_2$ in $V_4$ is a vector subspace of dimension $d_{\mathfrak m}\in\{0,1,2\}$. The images of $(a,b,c,d)\in Um(R^4)$ and $(d^{\prime},-c^{\prime},-b^{\prime},a^{\prime})\in Um(R^4)$ in $V_2$ are two nonzero vectors $v$ and $v^{\prime}$ (respectively), hence $d_{\mathfrak m}\ge 1$. We show that the assumption that $d_{\mathfrak m}=1$ leads to a contradiction. This assumption implies that there exists a nonzero scalar $\beta\in R/\mathfrak m$ such that $v=\beta v^{\prime}$. If $\zeta_{\mathfrak m}$ is the image of $\zeta$ in $V_4$, by denoting also by $\langle\,,\rangle$ its reduction modulo $\mathfrak m$, we compute $0=\langle \zeta_{\mathfrak m},v^{\prime}\rangle=\beta\langle \zeta_{\mathfrak m},v\rangle=\beta$, a contradiction to $\beta\neq 0$. Thus $d_{\mathfrak m}\neq 1$, hence $d_{\mathfrak m}=2$. As $d_{\mathfrak m}=2$ for each $\mathfrak m\in\textup{Max\,} R$, we have $M_2\cong R^2$ and $R^4/M_2$ is a projective $R$-module of rank $2$ by Lemma \ref{L3} below. Hence $M_2$ is a direct summand of $R^4$. 

Let 
$$\mathcal D^{\prime}:=R[X^{\prime},Y^{\prime},Z^{\prime},W^{\prime}]/(aX^{\prime}+bY^{\prime}+cZ^{\prime}+dW^{\prime},d^{\prime}X^{\prime}-c^{\prime}Y^{\prime}-b^{\prime}Z^{\prime}+a^{\prime}W^{\prime}).$$
As $\Lambda_A(\zeta)=aa^{\prime}+bb^{\prime}+cc^{\prime}+dd^{\prime}=1$, we have $\textup{Spec\,} R=\cup _{f\in
\{a,b,c,d\}}\textup{Spec\,} R_{ff^{\prime}}$ (if $f=a$ then $f^{\prime}=a^{\prime}$, if $f=b$ then $f^{\prime}=b^{\prime}$, etc.). For $f\in\{a,d\}$ (resp.\ $f\in\{b,c\}$), the $R_{ff^{\prime}}$-algebra $\mathcal D^{\prime}_{ff^{\prime}}$ is isomorphic to $R_{ff^{\prime}}[Y^{\prime},Z^{\prime}]$ (resp.\ $R_{ff^{\prime}}[X^{\prime},W^{\prime}]$). Let $J$ and $J^{\prime}$ be the ideals of $\mathcal D$ and $\mathcal D^{\prime}$ (respectively) generated by the images of $X^{\prime}$, $Y^{\prime}$, $Z^{\prime}$, $W^{\prime}$. We view $\mathcal D$ and $\mathcal D^\prime$ as augmented $R$-algebras (in the terminology of \cite{BCW}, Sect.\ 4, Cor.\ (4.3)) with the augmentations $\mathcal D\rightarrow R$ and $\mathcal D^{\prime}\rightarrow R$ which are the $R$-algebra homomorphisms having $J$ and $J^{\prime}$ (respectively) as their kernels. The $R$-modules $Q:=J/J^2$ and $Q^{\prime}=J^{\prime}/(J^{\prime})^2$ are identified via the Third Isomorphism Theorem to the following $R$-module quotient of ideals 
$$(X^{\prime},Y^{\prime},Z^{\prime},W^{\prime})/\bigl((aX^{\prime}+bY^{\prime}+cZ^{\prime}+dW^{\prime},d^{\prime}X^{\prime}-c^{\prime}Y^{\prime}-b^{\prime}Z^{\prime}+a^{\prime}W^{\prime})+(X^{\prime},Y^{\prime},Z^{\prime},W^{\prime})^2\bigr)$$ 
of $R[X^{\prime},Y^{\prime},Z^{\prime},W^{\prime}]$, thus they are isomorphic to $R^4/M_2$. From \cite{BCW}, Cor.\ (4.3) and the isomorphisms of localizations $\mathcal D_{a}\cong R_a[Y_1,Z_1]$ and $\mathcal D^{\prime}_{aa^{\prime}}\cong R_a[Y^{\prime},Z^{\prime}]$ that involve variables that are linear polynomials in $X_1,Y_1,Z_1,W_1$ and $X^{\prime},Y^{\prime},Z^{\prime},W^{\prime}$ (respectively) and their analogs with $a$ replaced by $b$, $c$, or $d$, it follows that the augmented $R$-algebras $\mathcal D$ and $\mathcal D^{\prime}$ are isomorphic to the symmetric $R$-algebras of $Q$ and $Q^{\prime}$ (respectively) endowed with their natural augmentations, and so they are isomorphic. So the $R$-algebra $\mathcal D$ is isomorphic to the symmetric $R$-algebra of $M_2$. 

As $M_2$ is a direct summand of $R^4$, there exists a short exact sequence $0\rightarrow R\rightarrow P^*\rightarrow Q^{\prime}\rightarrow 0$ of $R$-modules, so $Q\oplus R\cong Q^{\prime}\oplus R\cong P^*\cong P$. Thus $Q\oplus R^2\cong P\oplus R\cong R^4$. As the $R$-module $Q$ is stably free of rank $2$, it is self-dual by \cite{lam}, Ch.\ III, Sect.\ 6, Thm.\ 6.8.
\end{proof}

For reader's convenience we include the following lemma whose statement and proof are probably well-known, but for which we could not find a reference.

\begin{lemma}\label{L3}
Let $(i,j)\in [\mathbb N\cup\{0\}]^2$. Let $M_i$ be an $R$-submodule of $R^{i+j}$ generated by $i$ elements and such that for each $\mathfrak m\in\textup{Max\,} R$ the image of $M_i$ in $R^{i+j}/\mathfrak mR^{i+j}$ is a vector space of dimension $i$ over the field $R/\mathfrak m$. Then $M_i\cong R^i$ and $R^{i+j}/M_i$ is a projective $R$-module of rank $j$.
\end{lemma}

\begin{proof}
Let $l:R^i\rightarrow M_i$ be a surjective $R$-linear map. For $\mathfrak m\in\textup{Max\,} R$, let $l_{\mathfrak m}^{\perp}:R^j\rightarrow R^{i+j}$ be an $R$-linear map such that the rule $(x_1,\ldots,x_{i+j})\mapsto l(x_1,\ldots,x_i)+l_{\mathfrak m}^{\perp}(x_{i+1},\ldots,x_{i+j})$ defines an $R$-linear map $\ell_{\mathfrak m}:R^{i+j}\rightarrow R^{i+j}$ whose reduction modulo $\mathfrak m$ is an isomorphism. Let $f_{\mathfrak m}\in R\setminus\mathfrak m$ be the determinant of $\ell_{\mathfrak m}$. As the $R_{f_{\mathfrak m}}$-linear map $(\ell_{\mathfrak m})_{f_{\mathfrak m}}$ is an isomorphism, it follows that $l_{f_{\mathfrak m}}$ is an isomorphism and $(R^{i+j}/M_i)_{f_{\mathfrak m}}\cong R^j_{f_{\mathfrak m}}$. As $\textup{Spec\,} R=\cup_{\mathfrak m\in\textup{Max\,} R} \textup{Spec\,} R_{f_{\mathfrak m}}$, we conclude that $l$ is an isomorphism and $R^{i+j}/M_i$ is a projective $R$-module of rank $j$.
\end{proof}

\begin{corollary}\label{C5}
Assume there exist two ideals $\mathfrak i_1$ and $\mathfrak i_2$ of $R$ such that $\mathfrak i_1\cap\mathfrak i_2=\{0\}$ and $\det(A)\in\mathfrak i_2$. Then for $\mathcal C\in\{\mathcal D,\mathcal W\}$, each $R$-algebra homomorphism $\mathcal C\rightarrow R/\mathfrak i_1$ lifts to an $R$-algebra homomorphism $\mathcal C\rightarrow R$.
\end{corollary}

\begin{proof}
Let $h_{1,2}:\mathcal C\rightarrow R/(\mathfrak i_{1}+\mathfrak i_{2})$ be induced by an $R$-algebra homomorphism $h_{1}:\mathcal C\rightarrow R/\mathfrak i_1$. As $A$ modulo $\mathfrak i_2$ has zero determinant, $\mathcal C/\mathfrak i_{2}\mathcal C$ is the symmetric algebra of a projective $R/\mathfrak i_{2}$-module $Q_{2}$ of rank $2$ if $\mathcal C=\mathcal D$ and of rank $3$ if $\mathcal C=\mathcal W$ (see Theorem \ref{TH6}(1) and (8)). The $R$-algebra homomorphism $h_{1,2}$ is
uniquely determined by an $R$-linear map $l_{1,2}:Q_{2}\rightarrow
R/(\mathfrak i_{1}+\mathfrak i_{2})$. If $l_{2}:Q_{2}\rightarrow R/\mathfrak i_{2}$ is an $R$-linear
map that lifts $l_{1,2}$ and if $h_{2}:\mathcal C\rightarrow R/\mathfrak i_{2}$ is
the $R$-algebra homomorphism uniquely determined by $l_2$, then $h_{2}$
lifts $h_{1,2}$. As $\mathfrak i_{1}\cap \mathfrak i_{2}=\{0\}$, we have a pullback diagram
$$\xymatrix@R=10pt@C=16pt{
R \ar[r] \ar[d] & R/\mathfrak i_1 \ar[d]\\
R/\mathfrak i_2 \ar[r] & R/(\mathfrak i_1+\mathfrak i_2),\\}$$ 
so there exists a unique $R$-algebra homomorphism $\mathcal C\rightarrow R$ that lifts $h_{1}$ and $h_{2}$. 
\end{proof}

\begin{corollary}\label{C6}
Assume $R$ is a Hermite ring. With $A=\left[ 
\begin{array}{cc}
a & b \\ 
c & d\end{array}\right]\in Um\bigl(\mathbb M_2(R)\bigr)$, let $P$ be as in Equation (\ref{EQP}). Then the following properties hold.

\medskip
{\bf (1)} We have $P\cong R^3$ and the $R$-algebra $\mathcal U$ is a polynomial $R$-algebra in $3$ variables.

\smallskip
{\bf (2)} Assume $\det(A)=0$ and let the projective $R$-module $Q$ be as in Theorem \ref{TH6}(8). Then $Q\cong R^2$ and the $R$-algebra $\mathcal D$ is a polynomial $R$-algebra in $2$ variables.\end{corollary}

\begin{proof} 
The stably free $R$-module $P$ is free by \cite{WW}, Cor.\ 3.2. As $P\cong R^3$, part (1) follows from Theorem \ref{TH6}(1). 

If $\det(A)=0$, then we similarly argue that the stably free $R$-module $Q$ is free and that part (2) follows from Theorem \ref{TH6}(8).\end{proof}

\section{Proof of Theorem \ref{TH1}}\label{S4}

For $E\in\mathbb M_2(R)$ let $\textup{Tr}(E)$ be its trace, let $\textup{adj}(E)\in\mathbb M_2(R)$ be its adjugate, and let $\textup{Ker}_E$ be the kernel and $\textup{Im}_E$ be the image of the $R$-linear map $L_E:R^2\rightarrow R^2$ defined by $E$. Let $I_2:=\left[ 
\begin{array}{cc}
1 & 0 \\ 
0 & 1\end{array}\right]\in \mathbb M_2(R)$. 

We first prove the following general lemma.

\begin{lemma}\label{L4}
Let $G,H,E\in \mathbb M_2(R)$. Then the following properties hold.

\medskip \textbf{(1)} There exists a matrix $O\in\mathbb M_2(R)$ such that $H=G(I_2+\textup{adj}(G)O)$ iff $G$ and $H$ are congruent modulo $R\det(G)$. 

\smallskip \textbf{(2)} If $GE$ is unimodular, then $G$ and $E$ are unimodular. 

\smallskip \textbf{(3)} If $G$ is unimodular and $G$ and $GE$ are congruent modulo $R\det(G)$, then $GE$ is unimodular iff $E$ is unimodular.
\end{lemma}

\begin{proof}
As $G\textup{adj}(G)=\det(G)I_2$, for $O\in\mathbb M_2(R)$ we have $H=G+\det(G)O$ iff $H=G(I_2+\textup{adj}(G)O)$. So part (1) holds. 

The only nontrivial implication of parts (2) and (3) is the `if' of part (3). It suffices to show that the ideal $\mathfrak h$ of $R$ generated by the entries of $GE$ is not contained in any $\mathfrak m\in\textup{Max\,} R$. This holds if $\det(G)\in\mathfrak m$ as $G$ and $GE$ are congruent modulo $R\det(G)$. If $\det(G)\notin\mathfrak m$, then $G$ modulo $\mathfrak m$ is invertible, thus $GE$ modulo $\mathfrak m$ is nonzero as this is so for $E$ modulo $\mathfrak m$, so $\mathfrak h\not\subseteq\mathfrak m$.\end{proof}

We are now ready to prove Theorem \ref{TH1}. Let $A=\left[ 
\begin{array}{cc}
a & b \\ 
c & d\end{array}\right]\in Um\bigl(\mathbb M_2(R)\bigr)$ be as in Theorem \ref{TH1}. To $A$ and a quadruple $\upsilon=(x,y,z,w)\in R^4$, we associate matrices in $\mathbb M_2(R)$ as follows:

$$C_{\upsilon}:=\left[\begin{array}{cc}
-w & z \\ 
y & -x\end{array}\right],\;\;\;\;\;\;\;\;\;\;\;\;\;\;\;\;\;\;\;\;\;\;\;\;\;\;\;\;\;\;\;\;\;\;\;\;\;\;\;\;\;\;\;\;\;\;\;\;\;\;\;\;\;\;\;\;\;\;\;\;\;\;\;\;\;\;\;$$
$$D_{A,\upsilon}:=I_2+\textup{adj}(A)C_{\upsilon}=\left[ 
\begin{array}{cc}
1-by-dw & bx+dz \\ 
ay+cw & 1-ax-cz\end{array}\right],\;\;\;\;\;\;\;\;\;\;\;\;\;\;\;\;$$
$$B_{A,\upsilon}:=AD_{A,\upsilon}=A+\det(A)C_{\upsilon}=\left[ 
\begin{array}{cc}
a-\det(A)w & b+\det(A)z \\ 
c+\det(A)y & d-\det(A)x\end{array}\right].$$
We have $\textup{Tr}(D_{A,\upsilon})=2-ax-by-cz-dw$ and one computes
\begin{equation}\label{EQ1a}
\det(D_{A,\upsilon})=\Phi_A(x,y,z,w),
\end{equation}
\begin{equation}\label{EQ1b}
\;\;\;\;\;\;\;\;\;\;\;\det(B_{A,\upsilon})=\det(A)\Phi_A(x,y,z,w),
\end{equation}
\begin{equation}\label{EQ4}
\,1-\det(A)\det(C_{\upsilon})=\textup{Tr}(D_{A,\upsilon})-\det(D_{A,\upsilon}).
\end{equation}

To prove Theorem \ref{TH1}, we first remark that clearly $(2)\Rightarrow (1)$ and $(2)\Rightarrow (4)$. 

To prove that $(3)\Rightarrow (2)$, let the quadruple $\upsilon=(x,y,z,w)\in R^4$ be such that $ax+by+cz+dw=1$
and $xw-yz=0$. We claim that the matrix $C:=C_{\upsilon}$ satisfies the conditions in (2). We have $\Phi(x,y,z,w)=0$, $\det(C)=0$, $C\in Um\bigl(\mathbb M_2(R)\bigr)$ and for $B_{A,\upsilon}=A+\det(A)C$ we have $\det(B_{A,\upsilon})=0$ by Equation (\ref{EQ1b}). As $\textup{Tr}(D_{A,\upsilon})-\det(D_{A,\upsilon})=1$ by Equation (\ref{EQ4}), $D_{A,\upsilon}$ is unimodular, so $B_{A,\upsilon}=AD_{A,\upsilon}$ is unimodular by Lemma \ref{L4}(3). Hence $(3)\Rightarrow (2)$.

To show that $(4)\Rightarrow (3)$, let $C\in\mathbb M_2(R)$ be such that $\det(C)=\det\bigl(A+\det(A)C\bigr)=0$. There exists a unique $\upsilon=(x,y,z,w)\in R^4$ such that $C=C_{\upsilon}$; so $A+\det(A)C=B_{A,\upsilon}$. We have $\det(B_{A,\upsilon})=\det(C_{\upsilon})=0$. Thus $xw-zy=0$ and $\det(A)\Phi(x,y,z,w)=\det(B_{A,\upsilon})=0$ by Equation (\ref{EQ1b}). If $\det(A)\notin Z(R)$, then $\Phi(x,y,z,w)=0$; from this and $xw-yz=0$ it follows that $1-ax-by-cz-dw=0$. Hence $(4)\Rightarrow (3)$ if $\det(A)\notin Z(R)$. 

In general, $(3)$ holds iff there exists an $R$-algebra homomorphism $h:\mathcal D\rightarrow R$ by  Lemma \ref{L1}, Equation (\ref{EQ3.e}) applied to $S=R$. By replacing $R$ with a finitely generated $\mathbb{Z}$-subalgebra $S$
of $R$ such that $A,B_{A,\upsilon}\in Um\bigl(\mathbb{M}_{2}(S)\bigr)$ and $C=C_{\upsilon}\in\mathbb M_2(S)$, to show that $h$ exists we can assume that $R$ is noetherian. Thus $N(R)$ is nilpotent and there exists $j\in\mathbb N$ such that the set $\{\mathfrak{p}
_{1},\ldots ,\mathfrak{p}_{j}\}$ of minimal prime ideals of $R$ has $j$ elements. As the $R$-algebra $\mathcal D$ is smooth (see Theorem \ref{TH6}(3)) and $N(R)$ is nilpotent, each $R$-algebra homomorphism $\mathcal D\rightarrow R/N(R)$ lifts to an $R$-algebra homomorphism $\mathcal D\rightarrow R$ (e.g., see \cite{BLR}, Ch.\ 2, Sect.\ 2.2, Prop.\ 6). Thus, by replacing $R$ with $R/N(R)$, we can assume also that $N(R)=\cap
_{i=1}^{j}\mathfrak{p}_{i}=\{0\}$. If $\det(A)=0$, then the $R$-algebra $\mathcal D$ is isomorphic to a symmetric $R$-algebra by Theorem \ref{TH6}(8) and thus $h$ exists. So we can assume also that $\det(A)\neq 0$. Let $\det(A) _{i}:=\det(A)+\mathfrak{p}_{i}\in R/\mathfrak{p}_{i}$. As $\det(A) \neq 0$, there exists an index $i\in \{1,\ldots ,j\}$ such that $\det(A)\notin\mathfrak p_i$, i.e., $\det(A) _{i}\neq 0$. We can assume that the minimal prime
ideals are indexed such that there exists $j^{\prime }\in \{1,\ldots ,j\}$ for which $\det(A) _i\neq 0$ if $i\in\{1,\ldots,j^{\prime}\}$ and $\det(A) _i=0$ if $i\in\{j^{\prime}+1,\ldots,j\}$. If $\mathfrak i_{1}:=\cap
_{i=1}^{j^{\prime }}\mathfrak{p}_{i}$ and 
$\mathfrak i_{2}:=\cap _{i=j^{\prime }+1}^{j}\mathfrak{p}_{i}$, we have 
$\mathfrak i_{1}\cap \mathfrak i_{2}=\{0\}$ and $\det(A)\in \mathfrak i_2$. As $\det(A)+\mathfrak i_1\notin Z(R/\mathfrak i_{1})$, from the prior paragraph applied to $R/\mathfrak i_1$ and  Lemma \ref{L1}, Equation (\ref{EQ3.e}) applied to $S=R/\mathfrak i_1$ it follows that there exists an $R$-algebra homomorphism $h_{1}:\mathcal D\rightarrow R/\mathfrak i_{1}$. So there exists $h:\mathcal D\rightarrow R$ that lifts $h_1$ by Corollary \ref{C5}. Hence $(4)\Rightarrow (3)$. 

We conclude that $(1)\Leftarrow (2)\Leftrightarrow (3)\Leftrightarrow (4)$.

We prove that $(1)\Rightarrow (2)$. As $(2)\Leftrightarrow (3)$, as above we argue that it suffices to prove that $(1)\Rightarrow (2)$ when $R$ is noetherian and $N(R)=\{0\}$. Let the ideals $\mathfrak i_1$ and $\mathfrak i_2$ of $R$ be as above. Let $B\in Um\bigl(\mathbb M_2(R)\bigr)$ be congruent to $A$ modulo $R\det(A)$ and $\det(B)=0$. Let $\upsilon=(x,y,z,w)\in R^4$ be such that $B=B_{A,\upsilon}$ (see Lemma \ref{L4}(1)). With $C:=C_{\upsilon}$ and $D:=D_{A,\upsilon}$, as $B=AD\in Um\bigl(\mathbb M_2(R)\bigr)$ we have $D\in Um\bigl(\mathbb M_2(R)\bigr)$ (see Lemma \ref{L4}(2)). As $\det(A)+\mathfrak i_1\notin Z(R/\mathfrak i_1)$, from the identity $\det(B)=\det(A)\det(D)=0$ and Equation (\ref{EQ1b}) it follows that $\Phi(x,y,z,w)\in\mathfrak i_1$, hence there exists an $R$-algebra homomorphism $g_1:\mathcal W\rightarrow R/\mathfrak i_1$ that maps the elements $\bar{X},\bar{Y},\bar{Z},\bar{W}$ of $\mathcal W$ to $x+\mathfrak i_1,y+\mathfrak i_1,z+\mathfrak i_1,w+\mathfrak i_1$ (respectively). As $\mathfrak i_1\cap\mathfrak i_2=\{0\}$, let $g:\mathcal W\rightarrow R$ be an $R$-algebra homomorphism that lifts $g_{1}$ (see Corollary \ref{C5}). 

Let $\upsilon^{\prime}=(x^\prime,y^{\prime},z^{\prime},w^{\prime}):=g^4(\bar{X},\bar{Y},\bar{Z},\bar{Z})\in R^4$. For the matrices $C^{\prime}:=C_{\upsilon^{\prime}}$ and $D^{\prime}:=D_{A,\upsilon^{\prime}}$ we have (see Equation (\ref{EQ1a})) $\det(D^\prime)=\Phi(x^{\prime},y^{\prime},z^{\prime},w^{\prime})=0$ and $C^{\prime}$ and $C$ are congruent modulo $\mathfrak i_1$. Hence $D^{\prime}$ and $D$ are congruent modulo $\mathfrak i_1$. As $D$ is unimodular, it follows that the ideal $\mathfrak d^\prime$ of $R$ generated by the entries of $D^\prime$ satisfies $\mathfrak d^{\prime}+\mathfrak i_1=R$ and thus $\mathfrak d^{\prime}$ is not contained in any $\mathfrak m\in\textup{Max\,} R$ with $\det(A)\notin\mathfrak m$. As $\textup{Tr}(D^{\prime})-\det(D^{\prime})=1-\det(A)\det(C^{\prime})\in\mathfrak d^{\prime}$ by Equation (\ref{EQ4}), $\mathfrak d^{\prime}$ is not contained in any maximal ideal which does not contain $1-\det(A)\det(C^{\prime})$. The last two sentences imply that $\mathfrak d^{\prime}$ is not contained in any $\mathfrak m\in\textup{Max\,} R$, thus $\mathfrak d^{\prime}=R$, i.e., $D^{\prime}$ is unimodular. From Lemma \ref{L4}(3) it follows that $B^{\prime}:=AD^{\prime}$ is unimodular. 

By replacing the triple $(C,B,D)$ with $(C^{\prime},B^{\prime},D^{\prime})$, we can assume that $\det(D)=0$. As $D=I_2+\textup{adj}(A)C$ has zero determinant, it follows that $C\in Um\bigl(\mathbb M_2(R)\bigr)$. 

To complete the proof that $(1)\Rightarrow (2)$, it suffices to show that we can replace $C$ by a matrix $C_1\in Um\bigl(\mathbb M_2(R)\bigr)$ with $\det(C_1)=0$ and such that for $D_1:=I_2+\textup{adj}(A)C_1$ we have $\det(D_1)=0$ and $D_1\in Um\bigl(\mathbb M_2(R)\bigr)$: so $B_1:=AD_1$ is congruent to $A$ modulo $R\det(A)$ and unimodular by Lemma \ref{L4}(1) and (3) with $\det(B_1)=0$. As $\textup{Ker}_D$ and $\textup{Im}_D$ are projective $R$-modules of rank $1$ (see \cite{CPV1}, Lem.\ 3.1), the short exact sequence $0\rightarrow\textup{Ker}_D\rightarrow R^2\rightarrow\textup{Im}_D\rightarrow 0$ splits, i.e., it has a section $\sigma:\textup{Im}_D\rightarrow R^2$. Let $C_1\in\mathbb M_2(R)$ be the unique matrix such that $\textup{Ker}_{D}\subseteq\textup{Ker}_{C_1-C}$ and $\sigma(\textup{Im}_D)\subseteq\textup{Ker}_{C_1}$. As $\textup{Ker}_D$ is a direct summand of $R^2$ of rank $1$ and for $t\in\textup{Ker}_{D}$ we have $\textup{adj}(A)C_1(t)=\textup{adj}(A)C(t)=-t$, it follows first that $\textup{Ker}_D\subseteq\textup{Ker}_{D_1}$, second that
$\textup{Im}_{C_1}=C_1(\textup{Ker}_D)=C(\textup{Ker}_D)$ is a direct summand of $R^2$ of rank $1$ isomorphic to $\textup{Ker}_D$, and third that $\textup{Ker}_{C_1}=\sigma(\textup{Im}_D)$ is also a direct summand of $R^2$ of rank $1$. As $\textup{Ker}_D\subset\textup{Ker}_{D_1}$, $R^2=\textup{Ker}_D\oplus \sigma(\textup{Im}_D)$ and $\sigma(\textup{Im}_D)\subset\textup{Ker}_{C_1}$ we compute
$$\textup{Im}_{D_1}=D_1\bigl(\sigma(\textup{Im}_D)\bigr)=\{x+\textup{adj}(A)C_1(x)|x\in \sigma(\textup{Im}_D)\}=\sigma(\textup{Im}_D).$$
We conclude that $C_1,D_1\in Um\bigl(\mathbb M_2(R)\bigr)$ and $\det(C_1)=\det(D_1)=0$. Hence $(1)\Rightarrow (2)$, thus Theorem \ref{TH1} holds.

\section{Proof of Theorem \ref{TH2}}\label{S5}

Let $A\in Um\bigl(\mathbb M_2(R)\bigr)$. We include two proofs of Theorem \ref{TH2}(1). First, part (1) holds as clearly statement (4) of \cite{CPV1}, Thm.\ 4.3 implies statement (3) of Theorem \ref{TH1}. Second, recall the $R$-algebra homomorphism $\rho:\mathcal D\rightarrow\mathcal X$ of Equation (\ref{EQ3}); as $A$ is simply extendable there exists an $R$-algebra homomorphism $\mathcal X\rightarrow R$ by Lemma \ref{L2}, (1) implies (3), which composed with $\rho$ gives an $R$-algebra homomorphism $\mathcal D\rightarrow R$, hence statement (3) of Theorem \ref{TH1} holds by  Lemma \ref{L1}, Equation (\ref{EQ3.e}) and thus $A$ is determinant liftable. 

To prove part (2), assume that $A$ is extendable. Then $A$ modulo $R\det (A)$ is non-full by \cite{CPV1}, Prop.\ 5.1(2). Let $(\bar{l},\bar{m},\bar{o},\bar{q})\in [R/R\det (A)]^4$ be such that $A$ modulo $R\det (A)$
is $\left[ 
\begin{array}{c}
\bar{l} \\ 
\bar{m}\end{array}\right] \left[ 
\begin{array}{cc}
\bar{o} & \bar{q}\end{array}\right] $. If $(l,m,o,q)\in R^4$ lifts $(\bar{l},\bar{m},\bar{o},\bar{q})$, then $B:=\left[ 
\begin{array}{c}
l \\ 
m\end{array}\right] \left[ 
\begin{array}{cc}
o & q\end{array}\right] $ is congruent to $A$ modulo $R\det (A)$ and $\det (B)=0$, hence $A$ is weakly determinant liftable. Thus Theorem \ref{TH2} holds.

\begin{example}\label{EX3.5}
\normalfont
Let $R$ be a Dedekind domain which is not a PID. So $R$ is not a $\Pi_2$ ring by \cite{CPV1}, Thm.\ 1.7(4). Hence there exist a matrix $A\in Um\bigl(\mathbb M_2(R)\bigr)$ of zero determinant which is not (simply) extendable. As $A$ is (weakly) determinant liftable, it follows that the converses of Theorem \ref{TH3} do not hold. 
\end{example}

\section{A criterion for weakly determinant liftability}\label{S6}

As the $R$-algebra $\mathcal W$ of Section \ref{S2} is not smooth in general (see Theorem \ref{TH6}(5) and its proof), the analog of the equivalence $(1)\Leftrightarrow (3)$ in Theorem \ref{TH1} for weakly determinant liftability has the following more complex form. 

\begin{theorem}\label{TH9}
For $A=\left[ 
\begin{array}{cc}
a & b \\ 
c & d\end{array}\right]\in Um\bigl(\mathbb M_2(R)\bigr)$ let $\Phi=\Phi_A$ be as in Section \ref{S2}. Then the following properties hold.

\medskip \textbf{(1)} If $A$ is determinant liftable, then there exists $(x,y,z,w)\in R^4$ such that $\Phi(x,y,z,w)=0$.

\smallskip \textbf{(2)} If there exists $(x,y,z,w)\in R^4$ such that $\Phi(x,y,z,w)=0$, then $A$ is weakly determinant liftable.

\smallskip \textbf{(3)} If either $N(R)=\{0\}$ or $\det(A)\notin Z(R)$, then the converse of part (2) holds.

\smallskip \textbf{(4)} If $\det(A)\in Z(R)$ and $A$ is weakly determinant liftable, then there exists $(x,y,z,w)\in R^4$ such that $\Phi(x,y,z,w)\in N(R)$.
\end{theorem}

\begin{proof}
If $A$ is determinant liftable, then there exists $(x,y,z,w)\in R^4$ with $xw-yz=0$ and $ax+by+cz+dw=1$ by Theorem \ref{TH1}, so $\Phi(x,y,z,w)=0$. So part (1) holds. 

If $\upsilon=(x,y,z,w)\in R^4$ is as in part (2), then for $B:=B_{A,\upsilon}$ we have $\det(B)=0$ by Equation (\ref{EQ1b}). Thus, as $A$ and $B$ are congruent modulo $R\det(A)$, $A$ is weakly determinant liftable. So part (2) holds.

To prove part (3) let $B\in\mathbb M_2(R)$ be congruent to $A$ modulo $R\det(A)$ with $\det(B)=0$. Let $\upsilon=(x,y,z,w)\in R^4$ be such that $B=B_{A,\upsilon}=A\bigl(I_2+\textup{adj}(A)C_{\upsilon}\bigr)$ (see Lemma \ref{L4}(1)). Thus $\det(A)\Phi(x,y,z,w)=\det(B)=0$ by Equation (\ref{EQ1b}). If $\det(A)\notin Z(R)$, then from $\det(A)\Phi(x,y,z,w)=0$ it follows that $\Phi(x,y,z,w)=0$. So we can assume that $N(R)=\{0\}$. As in Section \ref{S4} we argue that we can assume that $R$ is noetherian; if the ideals $\mathfrak i_1$ and $\mathfrak i_2$ of $R$ are as in Section \ref{S4}, then $\det(A)\in\mathfrak i_2$ and $\det(A)+\mathfrak i_1\notin Z(R/\mathfrak i_1)$. As $\det(A)+\mathfrak i_1\notin Z(R/\mathfrak i_1)$ and $\det(A)\Phi(x,y,z,w)=0$, it follows that $\Phi(x,y,z,w)\in\mathfrak i_1$. So as in the first paragraph of the proof of the implication $(1)\Rightarrow (2)$ of Section \ref{S4} we argue based on Corollary \ref{C5} that there exists an $R$-algebra homomorphism $\mathcal W\rightarrow R$. From this and  Lemma \ref{L1}, Equation (\ref{EQ3.d}) it follows that there exists $(x,y,z,w)\in R^4$ with $\Phi(x,y,z,w)=0$. So part (3) holds.

Part (4) follows from part (3) applied to $R/N(R)$.\end{proof}

\begin{example}\label{EX4}
\normalfont
If $R$ is such that $N(R)=\{0\}$ and there exists $A\in Um\bigl(\mathbb M_2(R)\bigr)$ which is not determinant liftable but is weakly determinant liftable (see Example \ref{EX1}), then there exists $(x,y,z,w)\in R^4$ such that $\Phi(x,y,z,w)=0$ by Theorem \ref{TH9}(3). Hence the converse of Theorem \ref{TH9}(1) does not hold in general.
\end{example}

\begin{remark}\label{rem2}
\normalfont

If $\upsilon=(x,y,z,w)\in R^4$ is such that $\Phi(x,y,z,w)\neq 0=\Phi(x,y,z,w)^2$ and the matrix $B_{A,\upsilon}$ is not unimodular, i.e., the ideal $\mathfrak b$ generated by its entries is not $R$, then there exists $(x^{\prime},y^{\prime},z^{\prime},w^{\prime})\in \upsilon + [R\Phi(x,y,z,w)]^4$ such that $\Phi(x^{\prime},y^{\prime},z^{\prime},w^{\prime})=0$ iff the annihilator of $\Phi(x,y,z,w)$ in $R$ contains an element of $1-\mathfrak b$ (i.e., iff we have $\Phi(x,y,z,w)\in\Phi(x,y,z,w)\mathfrak b$).
\end{remark}

\section{Proofs of Theorems \ref{TH3} and \ref{TH4}}\label{S7}

Assume the ring $R$ is such that the simply extendable and determinant liftable properties on a matrix in $Um\bigl(\mathbb M_2(R)\bigr)$ are equivalent. As unimodular matrices in $\mathbb M_2(R)$ of zero determinant are determinant liftable, they are simply extendable, so $R$ is a $\Pi_2$ ring by definition. Thus the `if' part of Theorem \ref{TH3} holds.

Based on Theorem \ref{TH2}(1), to prove the `only if' part of Theorem \ref{TH3} it suffices to show that if $R$ is a $\Pi_2$ ring and if for $A=\left[ 
\begin{array}{cc}
a & b \\ 
c & d\end{array}\right]\in Um\bigl(\mathbb M_2(R)\bigr)$ there exists $B=\left[ 
\begin{array}{cc}
a_1 & b_1 \\ 
c_1 & d_1\end{array}\right]\in Um\bigl(\mathbb{M}_{2}(R)\bigr)$ congruent to $A$ modulo $R\det(A)$ and $\det(B)=0$, then $A$ is simply extendable. As $R$ is a $\Pi_2$ ring, $B$ is simply extendable. From this and \cite{CPV1}, Thm.\ 4.3, (3) applied to $B$, it follows that there exists $(e,f)\in Um(R^2)$ such that $(a_1e+c_1f,b_1e+d_1f)\in Um(R^2)$ and so $(a_1e+c_1f,b_1e+d_1f, ad-bc)\in Um(R^3)$. As $B-A\in\mathbb M_2(R\det(A))$, it follows that $(ae+cf,be+df, ad-bc)\in Um(R^3)$. Thus $A$ is simply extendable by \cite{CPV1}, Cor.\ 4.7(2). So Theorem \ref{TH3} holds. 

The `only if' of Theorem \ref{TH4} holds by Theorem \ref{TH2}(2). 

To prove the `if' part of Theorem \ref{TH4}, we assume that $A$ is weakly determinant liftable. Let $B\in\mathbb{M}_{2}(R)$ be congruent to $A$ modulo $R\det(A)$ and $\det(B)=0$. As $B$ is non-full by hypothesis, $A$ modulo $R\det(A)$ is non-full and thus $A$ is extendable by \cite{CPV1}, Prop.\ 5.1(2). So Theorem \ref{TH4} holds. 

\section{On $WJ_{2,1}$ and $J_{2,1}$ rings}\label{S8}

We first prove Theorem \ref{TH5}. Let $R$ be a $WJ_{2,1}$ ring. Let $A =\left[ 
\begin{array}{cc}
a & b \\ 
c & d\end{array}\right]\in Um\bigl(\mathbb M_2(R)\bigr)$. By taking $(\Psi,\Delta)=(1,0)$ in Definition \ref{def2}(1), it follows that there exists $(x,y,z,w)\in R^4$ such that $ax+by+cz+dw=1$ and $xw-yz=0$, and hence $A$ is determinant liftable by Theorem \ref{TH1}. So part (1) holds. 

We assume now that $R$ is also a Hermite ring. The `only if' of part (2) follows from the fact that each elementary divisor ring is an $SE_2$ ring (see \cite{CPV1}, Prop.\ 1.3) and hence a $\Pi_2$ ring. For the `if' of part (2), if $R$ is also a $\Pi_2$ ring, then from part (1) and Corollary \ref{C2} we get that $R$ is an $SE_2$ ring and thus an elementary divisor ring by \cite{CPV1}, Cor.\ 1.8. So part (2) holds. 

Each Hermite domain is a B\'{e}zout domain and hence a pre-Schreier domain (see \cite{CPV1}, Sect.\ 2) and a $\Pi_2$ domain (see \cite{CPV1}, paragraph after Thm.\ 1.4). Based on this, part (3) follows from part (2). Thus Theorem \ref{TH5} holds.

\begin{proposition}\label{P1}
A ring $R$ is a $J_{2,1}$ ring in the sense of Definition \ref{def2}(2) iff it is a $J_{2,1}$ ring in the sense of \cite{lor}, Def.\ 4.6.
\end{proposition}

\begin{proof}
Assume $R$ is a $J_{2,1}$ ring in the sense of Definition \ref{def2}(2). Let $(\alpha,\beta,\gamma,\delta)\in R^4$. As $R$ is a Hermite ring, there exist $e\in R$ and $(a,b,c,d)\in Um(R^4)$ such that $(\alpha,\beta,\gamma,\delta)=e(a,b,c,d)$. For $\psi\in R$, the equation $\alpha X+\beta Y+\gamma Z+\delta W=\psi$ has a solution $(x,y,z,w)\in R^4$ iff $\Psi\in Re$. Let $(\Psi,\Delta)\in Re\times R$. Let $f\in R$ be such that $\Psi=ef$. From Definition \ref{def2} applied to $(\Psi,\Delta)=(f,0)$ it follows that there exists $(x,y,z,w)\in R^4$ such that $ax+by+cz+dw=f$ and $xw-yz=\Delta$. So $\alpha x+\beta y+\gamma z+\delta w=ef=\Psi$ and $xw-yz=\Delta$. Thus the `only if' part holds. 

Assume $R$ is a $J_{2,1}$ ring in the sense of \cite{lor}, Def.\ 4.6. Clearly, $R$ is a $WJ_{2,1}$ ring. As $R$ is a Hermite ring by \cite{lor}, Prop.\ 4.11, the `if' part holds. 
\end{proof}

\section{Rings with universal (weakly) determinant liftability}\label{S9} 

Let $GL_2(R)$ be the group of units of $\mathbb M_2(R)$. For a matrix $E\in\mathbb M_2(R)$, let $[E]\in {GL}_2(R)\backslash \mathbb{M}_2(R)/GL_2(R)$ be its double coset. For a projective $R$-module $M$ of rank $1$, let $[M]\in\textup{Pic\,}(R)$ be its class.

\begin{proposition}\label{P2}
We consider the following statements on $R$.

\medskip \textbf{(1)} For each $a\in R$, the map of sets
$$\{B\in Um\bigl(\mathbb M_2(R)\bigr)|\det(B)=0\}\rightarrow \{\bar{B}\in Um\bigl(\mathbb M_2(R/Ra)\bigr)|\det(\bar{B})=0\},$$
defined by the reduction modulo $Ra$, is surjective.

\smallskip \textbf{(2)} For each $a\in R$, the map of sets of double coset
$$\{[B]| B\in Um\bigl(\mathbb M_2(R)\bigr),\ \det(B)=0\}\rightarrow \{[\bar{B}]|\bar B\in Um\bigl(\mathbb M_2(R/Ra)\bigr),\ \det(\bar{B})=0\},$$
defined by the reduction modulo $Ra$, is surjective.

\smallskip \textbf{(3)} For each $a\in R$, every projective $R/Ra$-module of rank $1$ generated by $2$ elements is isomorphic to the reduction modulo $Ra$ of a projective $R$-module of rank $1$ generated by $2$ elements.

\smallskip \textbf{(4)} Each matrix in $Um\bigl(\mathbb M_2(R)\bigr)$ is determinant liftable.

\medskip
Then $(1)\Rightarrow (2)\Leftrightarrow (3)$ and $(1)\Rightarrow (4)$. If $sr(R)\le 4$, then $(1)\Leftrightarrow (4)$.
\end{proposition}

\begin{proof}
For a pair $\mu:=(M,M^{\prime})$ of projective $R$-submodules of $R^2$ of rank $1$ and generated by $2$ elements such that we have a direct sum decomposition $R^2=M\oplus M^{\prime}$, let $E_{\mu}\in Um\bigl(\mathbb M_2(R)\bigr)$ be the projection on $M$ along $M^{\prime}$; so $\det(E_{\mu})=0$, $M$ and $M^{\prime}$ are dual to each other (i.e., $[M^{\prime}]=-[M]$, with $\textup{Pic\,}(R)$ viewed additively), and $\mu_{[M]}:=[E_{\mu}]$ depends only on $[M]$. Each projective $R$-module of rank $1$ generated by $2$ elements is isomorphic to such an $M$. For $G\in Um\bigl(\mathbb M_2(R)\bigr)$ with $\det(G)=0$, $\textup{Ker}_G$ and $\textup{Im}_G$ are projective $R$-module of rank $1$ generated by $2$ elements and the short exact $0\to\textup{Ker}_G\to R^2\to\textup{Im}_G\to 0$ has a section $\varphi:\textup{Im}_G\rightarrow R^2$ (see \cite{CPV1}, Lem.\ 3.1); if $\tau_G:=\bigl(\textup{Im}(\varphi),\textup{Ker}_G\bigr)$, then $[G]=[E_{\tau_G}]=\mu_{[\textup{Im}_G]}$. Thus 
$$\{[B]| B\in Um\bigl(\mathbb M_2(R)\bigr),\ \det(B)=0\}=\{\mu_{[M]}|M\oplus M^{\prime}=R^2,\; M\;\textup{has rank}\;1\}.$$
From this and its analog over $R/Ra$, it follows that $(2)\Leftrightarrow (3)$. 

Clearly, $(1)\Rightarrow (2)$.

For $(1)\Rightarrow (4)$, let $A\in Um\bigl(\mathbb M_2(R)\bigr)$. By applying (1) to $a=\det(A)$ and the reduction $\bar B$ of $A$ modulo $Ra$, it follows that there exists $B\in Um\bigl(\mathbb M_2(R)\bigr)$ congruent to $A$ modulo $R\det(A)$ and $\det(B)=0$, so $A$ is determinant liftable. So $(1)\Rightarrow (4)$.

Assume $sr(R)\le 4$. To prove $(4)\Rightarrow (1)$, let $a\in R$. Let $\bar B\in Um\bigl(\mathbb M_2(R/Ra)\bigr)$ with $\det(\bar{B})=0$. Let $C\in Um\bigl(\mathbb M_2(R)\bigr)$ be such that its reduction modulo $Ra$ is $\bar B$ by \cite{CPV1}, Prop.\ 2.4(1); we have $\det(C)\in Ra$. As $C$ is determinant liftable, there exists $B\in Um\bigl(\mathbb M_2(R)\bigr)$ with $\det(B)=0$ and congruent to $C$ modulo $R\det(C)$ and hence also modulo $Ra$; so the map of statement (1) is surjective, hence $(4)\Rightarrow (1)$. 
\end{proof}

\begin{example}\label{EX5}
\normalfont
If $R$ is an integral domain of dimension 1, then each matrix $A\in Um\bigl(\mathbb M_2(R)\bigr)$ is determinant liftable. To check this we can assume that $\det(A)\neq 0$ and this case follows from \cite{CPV1}, Thm.\ 1.7(1) and Theorem \ref{TH2}(1). Recall that $sr(R)\le 2$ (cf.\ \cite{CPV1}, Sect.\ 1). So parts (1) to (4) of Proposition \ref{P2} hold.
\end{example}

\begin{proposition}\label{P3}
We consider the following two statements on $R$.

\medskip \textbf{(1)} For each $a\in R$, $Um\bigl(\mathbb M_2(R/Ra)\bigr)$ is contained in the image of the modulo $Ra$ reduction map
$\{B\in \mathbb M_2(R)|\det(B)=0\}\rightarrow \{\bar{B}\in \mathbb M_2(R/Ra)|\det(\bar{B})=0\}$.

\smallskip \textbf{(2)} Each matrix in $Um\bigl(\mathbb M_2(R)\bigr)$ is weakly determinant liftable.

\medskip
Then $(1)\Rightarrow (2)$, and the converse holds if $sr(R)\le 4$.
\end{proposition}

\begin{proof} It is the same as the last two paragraphs of the proof of Proposition \ref{P2}, with determinant and $B\in Um\bigl(\mathbb M_2(R)\bigr)$ replaced by weakly determinant and $B\in\mathbb M_2(R)$ (respectively).
\end{proof}

\section{A criterion for determinant liftability via completions}\label{S10}

The following proposition is probably well-known.

\begin{proposition}\label{P4} Let $A\in Um\bigl(\mathbb{M}_{2}(R)\bigr)$, let $t\in R$ be such that $\det
(A)\in Rt$ and let $\hat{R}$ be the $t$-adic completion of $R$. Then there
exists $B\in Um\bigl(\mathbb{M}_{2}(\hat{R})\bigr)$ whose reduction modulo $\textup{Ker}(\hat{R}\rightarrow R/Rt)$ is the reduction of $A$ modulo $Rt$ and $\det(B)=0$.
\end{proposition}

\begin{proof}
Let $B_0:=A$. By induction on $n\in \mathbb{N}$, we show that there exists $B_{n}\in 
\mathbb{M}_{2}(R)$ congruent to $B_{n-1}$ modulo $Rt^{2^{n-1}}$ and $\det (B_{n})\in
Rt^{2^{n}}$. For $n=1$, let $s\in R$ be such that $\det (A)=st$. With $A=\left[ 
\begin{array}{cc}
a & b \\ 
c & d\end{array}\right] $, for $B_{1}:=A+t\left[ 
\begin{array}{cc}
x & y \\ 
z & w\end{array}\right] \in \mathbb{M}_{2}(R)$, $\det (B_{1})$ is congruent modulo $Rt^{2}$
to $(dx+cy+bz+aw+s)t$. As $A\in Um\bigl(\mathbb{M}_{2}(R)\bigr)$, we can choose $(x,y,z,w)\in
R^{4}$ such that $dx+cy+bz+aw=-s$, hence $\det (B_{1})\in Rt^{2}$. The passage
from $n$ to $n+1$ follows from the case $n=1$ applied to $(B_{n},Rt^{2^{n+1}})$ instead of $(A,Rt^{2})$. This completes the induction. Let $B\in \mathbb{M}_{2}(\hat{R})$ be the limit of the sequence $(B_n)_{n\ge 1}$. Clearly, $\det(B)=0$. As $\textup{Ker}(\hat{R}\rightarrow R/Rt)$ is contained in the Jacobson radical of $\hat{R}$, we have $B\in Um\bigl(\mathbb M_2(\hat{R})\bigr)$. 
\end{proof}

Proposition \ref{P4} also follows from the smoothness part of Theorem \ref{TH6}(3) via a standard limit lifting argument. Proposition \ref{P4} gives directly the following result.

\begin{corollary}\label{C7}
Let $A\in Um\bigl(\mathbb{M}_{2}(R)\bigr)$. If $R$ is complete in the $\det (A)$-adic topology, then $A$ is determinant liftable.
\end{corollary}

\section{Applications of Picard groups}\label{S11}

To refine Theorem \ref{TH6}(7) in Theorem \ref{TH8}, we first recall that for a commutative $R$-algebra $S$ we have functorial (cocycle) isomorphisms (part of Hilbert's Theorem 90 for $\textup{Spec\,} R$)
\begin{equation}\label{EQ5}
\textup{Pic\,}(S)\cong H^1(\textup{Spec\,} S,\mathcal O_{\textup{Spec\,} S}^{\ast})\cong H^1(\textup{Spec\,} S,\mathbb{G}_{m,S})
\end{equation}
(e.g., see \cite{M}, Ch.\ III, Prop.\ 4.9), where $\mathcal O_{\textup{Spec\,} S}$ is the structure ringed sheaf on $\textup{Spec\,} S$, $\mathcal O_{\textup{Spec\,} S}^{\ast}$ is its set subsheaf of units and $H^1(\textup{Spec\,} S,\mathbb{G}_{m,S})$ is the group of equivalence classes of $\textup{Spec\,} S$-torsors under $\mathbb{G}_{m,S}$.

The following result which was first proved in \cite{CPV1}, (paragraph after) Thm.\ 1.4 is also a consequence of Theorem \ref{TH6}(8).

\begin{corollary}\label{C8}
If $\textup{Pic\,}(R)$ is trivial then $R$ is a $\Pi_2$ ring.
\end{corollary}

\begin{proof} We have to show each $A\in Um\bigl(\mathbb M_2(R)\bigr)$ with $\det(A)=0$ is simply extendable. Considering an $R$-algebra homomorphism $\varrho:\mathcal D\rightarrow R$ by Theorem \ref{TH6}(8), it defines a morphism $\textup{Spec\,} R\rightarrow\textup{Spec\,}\mathcal D$ of schemes and we pullback via it the $\textup{Spec\,}\mathcal D$-torsor under $\mathbb G_{m,\mathcal D}$ of Theorem \ref{TH6}(7). The resulting $\textup{Spec\,} R$-torsor under $\mathbb G_{m,R}$ is trivial as so is $\textup{Pic\,}(R)$ by Isomorphisms (\ref{EQ5}) and it is defined by a suitable action of $\mathbb G_{m,R}$ on $\textup{Spec\,} R\otimes_{\varrho,\mathcal D,\rho} \mathcal X\rightarrow\textup{Spec\,} R$. Thus there exists an $R$-algebra homomorphism $R\otimes_{\varrho,\mathcal D,\rho}\mathcal X\rightarrow R$ and therefore there exists also an $R$-algebra homomorphism $\mathcal X\rightarrow R$. From this and Lemma \ref{L2} we get that $A$ is simply extendable.
\end{proof}

For a commutative $R$-algebra $S$ defined by a homomorphism $\textup{Im}ath:R\rightarrow S$, let $\textup{Pic\,}(\textup{Im}ath):\textup{Pic\,}(R)\rightarrow\textup{Pic\,}(S)$ be the functorial homomorphism: for a projective $R$-module $M$ of rank $1$ we have $\textup{Pic\,}(\textup{Im}ath)([M]):=[S\otimes_R M]$. If $\jmath:S\rightarrow R$ is an $R$-algebra homomorphism, then as $\jmath\circ\textup{Im}ath$ is the identity automorphism of $R$, $\textup{Im}ath$ and $\textup{Pic\,}(\textup{Im}ath)$ are injective and $\jmath$ and $\textup{Pic\,}(\jmath)$ are surjective. In particular, if $S$ is a symmetric $R$-algebra then $\textup{Pic\,}(\textup{Im}ath)$ is injective.

For $A\in Um\bigl(\mathbb{M}_2(R)\bigr)$ and the $R$-algebra $\mathcal D=\mathcal D_A$ of Section \ref{S2}, let 
$$\iota_A:\textup{Pic\,}(R)\rightarrow\textup{Pic\,}(\mathcal D)$$ 
be the functorial homomorphism. If $A\in Um\bigl(\mathbb{M}_2(R)\bigr)$ is determinant liftable, then the $R$-algebra $\mathcal D$ has a retraction by Theorem \ref{TH1} and  Lemma \ref{L1}, Equation (\ref{EQ3.e}) and $\iota_A$ is injective. The $\textup{Spec\,}\mathcal D$-torsor under $\mathbb G_{m,\mathcal D}$ of Theorem \ref{TH6}(7) corresponds, under the isomorphism $\textup{Pic\,}(\mathcal D)\cong H^1(\textup{Spec\,}\mathcal D,\mathbb{G}_{m,\mathcal D})$, to a class $[\mathcal P]\in\textup{Pic\,}(\mathcal D)$, where $\mathcal P=\mathcal P_A$ is a projective $\mathcal D$-module of rank $1$; one would like to describe it and determine when it belongs to $\textup{Im}(\iota_A)$. 

We include a proof of the following well-known result\footnote{For instance, Lemma \ref{L5}(2) is stated in \cite{GH}, Sect.\ 0. Also, if $\mathfrak i$ is a nilradical ideal of $R$, then the pair $(R,\mathfrak i)$ is a henselian pair (e.g., see \url{https://stacks.math.columbia.edu/tag/09XD}) and hence Lemma \ref{L5}(2) is a particular case of \cite{GR}, Cor.\ 5.4.42 applied to $t=1$.\label{foo2}} as below it is essential.

\begin{lemma}\label{L5}
For an ideal $\mathfrak i$ of $R$ let $\pi:R\rightarrow R/\mathfrak i=:S$ be the quotient homomorphism. The following properties hold for $\iota:=\textup{Pic\,}(\pi):\textup{Pic\,}(R)\rightarrow\textup{Pic\,}(S)$.

\medskip
{\bf (1)} If $\mathfrak i$ is contained in the Jacobson radical of $R$, then $\iota$ is injective.

\smallskip
{\bf (2)} If $\mathfrak i\subset N(R)$, then $\iota$ is an isomorphism.
\end{lemma}

\begin{proof}
For part (1), let $M$ be a projective $R$-module of rank $1$ such that there exists an $S$-linear isomorphism $\bar\ell:S\rightarrow M/\mathfrak iM$. Let $m\in M$ be such that $\bar\ell(1+\mathfrak i)=m+\mathfrak iM$. Let $\ell:R\rightarrow M$ be the $R$-linear map such that $\ell(1)=m$. Then $M=\mathfrak iM+\ell(R)$. Thus $M=\ell(R)$ by Nakayama's Lemma (see \cite{bou}, Sect.\ 9, Subsect.\ 3, Thm.\ 2). So $\ell$ is surjective. As $M$ is a projective $R$-module, $\ell$ has an $R$-linear section, thus $R\cong M\oplus\textup{Ker}(\ell)$. From this and the fact that $M$ has rank $1$ it follows that $\textup{Ker}(\ell)$ is the zero $R$-module, so $\ell:R\rightarrow M$ is an isomorphism. Thus $\iota$ is injective.

Based on part (1), to prove part (2) it suffices to show that $\iota$ is surjective. We endow the set $\mathcal P_{\textup{f}}(\mathfrak i)$ of finite subsets $\Gamma$ of $\mathfrak i$ with the inclusion relation. For $\Gamma\in\mathcal P_{\textup{f}}(\mathfrak i)$, let $\mathfrak i_{\Gamma}$ be the ideal of $R$ generated by $\Gamma$. We get a direct system of $R$-algebras indexed by $\Gamma\in\mathcal P_{\textup{f}}(\mathfrak i)$ whose transition $R$-algebra homomorphisms for inclusions $\Gamma\subset\nabla$ with $\nabla\in\mathcal P_{\textup{f}}(\mathfrak i)$ are the natural surjections $R/\mathfrak i_{\Gamma}\rightarrow R/\mathfrak i_{\nabla}$. We have an $R$-algebra isomorphism $\varinjlim \{R/\mathfrak i_{\gamma}\}\rightarrow S$ and hence an isomorphism $\varinjlim \{\textup{Pic\,}(R/\mathfrak i_{\gamma})\}\rightarrow \textup{Pic\,}(S)$ by \cite{GH}, Thm.\ 1.3. Thus, by replacing $\mathfrak i$ with $\mathfrak i_{\Gamma}$s, we can assume that $\mathfrak i$ is finitely generated. So $\mathfrak i$ is nilpotent. For a projective $S$-module $\bar M$ of rank $1$, let $M$ be a projective $R$-module such that we have a surjective $R$-linear map $M\rightarrow\bar M$ whose kernel is $\mathfrak iM$ by \cite{bou}, Sect.\ 9, Subsect.\ 5, Prop.\ 11. Thus we have an $S$-linear isomorphism $M/\mathfrak iM\cong\bar M$. As $\bar M$ is finitely generated and $\mathfrak i$ is nilpotent, a similar argument based on Nakayama's Lemma gives that $M$ is finitely generated. As $\bar M$ has rank $1$ and $\textup{Spec\,} S\rightarrow\textup{Spec\,} R$ is a homeomorphism, $M$ has also rank $1$. As $\iota([M])=\bar M$, $\iota$ is surjective and part (2) holds.\end{proof}

For seminormal rings we refer to \cite{tr} and \cite{GH}. If $N(R)=\{0\}$ and we view $R$ as a subring of its total quotient ring $Q(R):=[R\setminus Z(R)]^{-1}R$, then (see \cite{GH}, Thm.\ 1.1) $R$ is seminormal iff it contains each $f\in Q(R)$ which is a root of a monic polynomial in $R[X]$ and for which there exists $n\in\mathbb N$ such that $f^r\in R$ for all integers $r\ge n$. A normal domain, i.e., an integral domain integrally closed in its field of fractions, is a seminormal domain. The next theorem only puts together several known results.

\begin{theorem}\label{TH7}
Let $S$ be a polynomial $R$-algebra in $n\in\mathbb N$ variables. The following properties hold for the injective functorial homomorphism $\iota:\textup{Pic\,}(R)\rightarrow\textup{Pic\,}(S)$.

\medskip
{\bf (1)} If $R/N(R)$ is not seminormal, then $\iota$ is not surjective.

\smallskip
{\bf (2)} Assume that for each $\mathfrak m\in\textup{Spec\,} R$, the local ring $R_{\mathfrak m}/N(R_{\mathfrak m})\cong \bigl(R/N(R)\bigr)_{\mathfrak m}$ is a seminormal domain. Then $\iota$ is an isomorphism

\smallskip
{\bf (3)} If $R$ is a Hermite ring, then $\iota$ is an isomorphism.
\end{theorem}

\begin{proof}
Based on Lemma \ref{L5}, we can assume that $N(R)=\{0\}$. 

See \cite{GH}, Thm.\ 1.5 for part (1). 

To prove parts (2) and (3), by replacing $R$ with $R_{\mathfrak m}$ with $\mathfrak m\in\textup{Max\,}(R)$, we can assume that $R$ is a local ring by \cite{qui}, Thm.\ 1.1. Thus part (2) holds by \cite{GH}, Thm.\ 1.6. For part (3), the local Hermite ring $R$ is a valuation ring by \cite{jen}, Thms.\ 1 and 2 and thus a valuation domain (this is stated in \cite{kap}, Sect.\ 10).\footnote{Recall from \cite{kap}, Sect.\ 10, Def.\ that a ring $R$ is called a valuation ring if for each $(a,b)\in R^2$, either $Ra\subset Rb$ or $Rb\subset Ra$, equivalently, if the ideals of $R$ are totally ordered by set inclusion. Valuation rings $S$ with $N(S)=\{0\}$ are integral domains. This is so as the assumption that there exists $(a,b)\in (S\setminus\{0\})^2$ such that $ab=0$ implies first that the nilpotent ideal $Sa\cap Sb$ is $\{0\}$ and second that the finitely generated ideal $Sa+Sb\cong Sa\oplus Sb$ is not principal, a contraction to $S$ being a valuation ring.\label{foo3}} Thus $R$ is a local normal domain and hence part (3) holds by part (2).\end{proof}

\begin{theorem}\label{TH8}
Let $A=\left[ 
\begin{array}{cc}
a & b \\ 
c & d\end{array}\right]\in Um\bigl(\mathbb M_2(R)\bigr)$. If $\det(A)\in N(R)$, then the following properties hold.

\medskip
{\bf (1)} For $f\in\{a,b,c,d\}$, the morphism $\textup{Spec\,}\mathcal X_f\rightarrow\textup{Spec\,}\mathcal D_f$ with the action of $\mathbb G_{m,\mathcal D_f}$ on it induced by the action of $\mathbb{G}_{m,\mathcal D}$ on $\textup{Spec\,} \rho:\textup{Spec\,}\mathcal X\rightarrow\textup{Spec\,}\mathcal D$ of Theorem \ref{TH6}(6) is a trivial $\textup{Spec\,}\mathcal D_f$-torsor under $\mathbb{G}_{m,\mathcal D_f}$.

\smallskip
{\bf (2)} The $\textup{Spec\,}\mathcal D$-torsor under $\mathbb{G}_{m,\mathcal D}$ of Theorem \ref{TH6}(7) is the pullback of a $\textup{Spec\,} R$-torsor under $\mathbb{G}_{m,R}$, i.e., we have $[\mathcal P]\in\textup{Im}(\iota_A)$.\end{theorem} 

\begin{proof} 
Based on Lemma \ref{L5}(2) and Isomorphisms (\ref{EQ5}), to prove parts (1) and (2) we can assume that $N(R)=\{0\}$; so $\det(A)=0$. For part (1) we only consider the case $f=a$ as the other three cases are similar. We have three $R_a$-algebra isomorphisms $\mathcal D_a\cong R_a[Y_1,Z_1]$ where $(Y_1,Z_1):=(aY+cW,Z+a^{-1}bW)$ (see proof of Theorem \ref{TH6}(8)) and (as $d=a^{-1}bc$ in $R_a$)
$$\mathcal X_a\cong R_a[X,Y,Z,W]/\bigl(1-(aX+cY)(W+a^{-1}bZ)\bigr)\cong R_a[X_2,Y,Z,W_2]/(1-X_2W_2)$$ 
for $(X_2,W_2):=(aX+cY,W+a^{-1}bZ)$. The action of $\mathbb{G}_{m,\mathcal D_a}$ on $\textup{Spec\,}\mathcal X_a$ defined by the $\mathcal D_a$-algebra homomorphism $\mathcal X_a\rightarrow\mathcal D_a[T,T^{-1}]\otimes_{\mathcal D_a,\rho_a} \mathcal X_a\cong\mathcal X_a[T,T^{-1}]$ that maps $\bar X$, $\bar Y$, $\bar Z$ and $\bar W$ to $T\bar X$, $T\bar Y$, $T^{-1}\bar Z$, $T^{-1}\bar W$ (see Section \ref{S2}), maps $\bar X_2$ to $T\bar X_2$ and $\bar W_2$ to $T^{-1}\bar W_2$, i.e., the substitution $(X_2,W_2)$ is compatible with the action. As $\rho_a(\bar Y_1)=\bar Z(a\bar X+c\bar Y)$ and $\rho_a(\bar Z_1)=\bar Y(\bar W+a^{-1}b\bar Z)$, under the three $R_a$-algebra isomorphisms, $\rho_a:\mathcal D_a\rightarrow\mathcal X_a$ gets identified to the $R_a$-algebra homomorphism $\varrho_a:R_a[Y_1,Z_1]\rightarrow R_a[X_2,Y,Z,W_2]/(1-X_2W_2)$ defined by $\varrho_a(Y_1):=\bar Z\bar X_2$ and $\varrho_a(Z_1):=\bar Y\bar W_2$. Clearly, $R_a[X_2,Y,Z,W_2]/(1-X_2W_2)\cong R_a[X_2,ZX_2,YW_2,W_2]/(1-X_2W_2)$. From either this and the compatibility of $(X_2,W_2)$ with the action or simply the existence of the $R_a[Y_1,Z_1]$-algebra homomorphism (retraction of $\varrho_a$) $R_a[X_2,Y,Z,W_2]/(1-X_2W_2)\rightarrow R_a[Y_1,Z_1]$ that maps $\bar X_2, \bar Y,\bar Z$ and $\bar W_2$ to $1,Z_1,Y_1$ and $1$ (respectively), it follows that part (1) holds.

For each $(f,g)\in\{a,b,c,d\}^2$, the functorial homomorphism $U(R_{fg})\rightarrow U(\mathcal D_{fg})$ is an isomorphism as $N(R_{fg})=\{0\}$ and $\mathcal D_{fg}$ is a polynomial $R_{fg}$-algebra by Theorem \ref{TH6}(8). For a projective $\mathcal D$-module $M$ of rank $1$ equipped with a $\mathcal D_f$-linear isomorphism $\lambda_f:M_f\cong\mathcal D_f$ for each $f\in\{a,b,c,d\}$, for every $(f,g)\in\{a,b,c,d\}^2$, the $\mathcal D_{fg}$-linear isomorphism $(\lambda_f)_g\circ (\lambda_g)_g^{-1}:\mathcal D_{fg}\rightarrow\mathcal D_{fg}$ is the multiplication by a unit $u_{f,g}\in U(R_{fg})$. Thus the cocyle $(u_{f,g})_{(f,g)\in\{a,b,c,d\}^2}$ that defines the class $[M]\in\textup{Pic\,}(\mathcal D)$ (see Isomorphisms (\ref{EQ5})), defines also a class in $\textup{Pic\,}(R)$ whose image under the functorial homomorphism $\textup{Pic\,}(R)\rightarrow\textup{Pic\,}(\mathcal D)$ is $[M]$. As the $\textup{Spec\,}\mathcal D$-torsor $\textup{Spec\,}\mathcal X\rightarrow\mathcal D$ under $\mathbb{G}_{m,\mathcal D}$ is defined via Isomorphisms (\ref{EQ5}) applied to $\mathcal D$ by such a class $[M]\in \textup{Pic\,}(\mathcal D)$ by part (1), part (2) holds.\end{proof} 

\begin{remark}
\normalfont

{\bf (1)} If $\det(A)\notin N(R)$, then one can check based on the proof of Theorem \ref{TH6}(8) that $\mathcal D$ is not a symmetric $R$-algebra.

\smallskip
{\bf (2)} Assume $\det(A)=0$. Let $\mathcal P_R=\mathcal P_{R,A}$ be a projective $R$-module of rank $1$ such that $[\mathcal P]=\iota_A([\mathcal P_R])$ by Theorem \ref{TH8}(2). One would like to describe all relations between the projective $R$-modules $P$ and $Q$ of Theorem \ref{TH6}(8) and the projective $R$-modules $\textup{Im}_A$ and $\mathcal P_R$ of rank $1$; recall from \cite{CPV1}, Lem.\ 3.1 that the image $\textup{Im}_A$ and the kernel $\textup{Ker}_A$ of the $R$-linear map $L_A:R^2\rightarrow R^2$ are projective $R$-modules of rank $1$ dual to each other.
\end{remark} 

\noindent \textbf{Acknowledgement.} The third author would like to thank Ofer Gabber for sharing the henselian pair part of Footnote \ref{foo2} and the argument of Footnote \ref{foo3} and SUNY Binghamton for good working conditions. The authors thank the referees for many valuable comments and suggestions.

\hbox{} \hbox{Grigore C\u{a}lug\u{a}reanu\;\;\;E-mail: calu@math.ubbcluj.ro}
\hbox{Address: Department of Mathematics, Babe\c{s}-Bolyai
University,} 
\hbox{1 Mihail Kogălniceanu Street, Cluj-Napoca 400084, Romania.}

\hbox{} \hbox{Horia F.\ Pop\;\;\;E-mail: horia.pop@ubbcluj.ro} 
\hbox{Address:
Department of Computer Science, Babe\c{s}-Bolyai
University,} 
\hbox{1 Mihail Kogălniceanu Street, Cluj-Napoca 400084, Romania.}

\hbox{} \hbox{Adrian Vasiu,\;\;\;E-mail: avasiu@binghamton.edu} 
\hbox{Address:
Department of Mathematics and Statistics, Binghamton University,} 
\hbox{P.\ O.\ Box
6000, Binghamton, New York 13902-6000, U.S.A.}

\end{document}